%
%
%   Franc Forstneric
%
%   THE HOMOTOPY PRINCIPLE IN COMPLEX ANALYSIS 
%
%   
%   For the Volume dedicated to Professor Robert Greene
%
%   January 8, 2003.
%
%   
%
%
\scrollmode
\magnification=\magstep1
\parskip=\smallskipamount

\hoffset=1cm \hsize=12cm

\def\demo#1:{\par\noindent\it{#1}. \rm}
\def\ni{\noindent}               % noindent
     % new indented par in horizontal mode
 % new unindented par in horizontal mode
\def\ll{\leftline}
\def\cl{\centerline}

%
%   \beginsection
%
\outer\def\beginsection#1\par{\bigskip
  \message{#1}\leftline{\bf\&#1}
  \nobreak\smallskip\vskip-\parskip\noindent}

%
%  \proclaim
%
\outer\def\proclaim#1:#2\par{\medbreak\vskip-\parskip
    \noindent{\bf#1.\enspace}{\sl#2}
  \ifdim\lastskip<\medskipamount \removelastskip\penalty55\medskip\fi}

\def\endpr{\hfill $\spadesuit$ \medskip}

%
%  Special letters: real, complex numbers etc.
%

\def\R{{\rm I\kern-0.2em R\kern0.2em \kern-0.2em}}
\def\N{{\rm I\kern-0.2em N\kern0.2em \kern-0.2em}}
\def\P{{\rm I\kern-0.2em P\kern0.2em \kern-0.2em}}
\def\B{{\rm I\kern-0.2em B\kern0.2em \kern-0.2em}}
\def\C{{\rm C\kern-.4em {\vrule height1.4ex width.08em depth-.04ex}\;}}
\def\CP{\C\P}
\def\RP{\R\P}

%
%
%  Roman capital characters
%
%
               % roman D in math mode

               % roman T
               % roman N in math

%
%
% Caligraphic capital characters
%
%

\def\cC{{\cal C}}

\def\cF{{\cal F}}
\def\cG{{\cal G}}

\def\cM{{\cal M}}

\def\cO{{\cal O}}

\def\cR{{\cal R}}

\def\cU{{\cal U}}

%
%
%  Small Greek letters in Math mode
%
\def\a{\alpha}
\def\b{\beta}
\def\g{\gamma}
\def\d{\delta}
\def\e{\epsilon}
\def\z{\zeta}

\def\l{\lambda}

\def\S{{\Sigma}}

%\def\O{\Omega}

%
%
%  Miscellaneous
%
%
\def\bar{\overline}              % conjugate
\def\bs{\backslash}              % backslash
             % disc
        % closed disc
                % partial derivative
\def\dibar{\bar\partial}         % di-bar derivative

%
% Abbreviations
%
\def\dim{{\rm dim}\,}                    % dimension
\def\holo{holomorphic}                   % holomorphic
                  % automorphism
                  % homomorphism
               % analytic subset
                % homeomorphism
\def\nbd{neighborhood}                   % neighborhood
                   % pseudoconvex
\def\spsc{strongly\ pseudoconvex}        % strongly psc
\def\ra{real-analytic}                   % real-analytic
               % plurisubharmonic
\def\spsh{strongly\ plurisubharmonic}
\def\tr{totally real}                    % totally real
             % polynomially convex
          % holo convex
           % holomorphic function
                % relatively compact subset
                  % support
             % C^n equivalent

\def\iff{if and only if}

\def\phe{proper holomorphic embedding}

\def\hra{\hookrightarrow}

\def\hvb{holomorphic vector bundle}

\def\Hom{{\rm Hom}}
                         % roman Aut in math mode
              % sypmlectic automorphisms

\def\begin{\ll{}\vskip 10mm \nopagenumbers}  % beginning of the paper
\def\pn{\footline={\hss\tenrm\folio\hss}}   % pagenumbers at bottom

\def\ii#1{\itemitem{#1}}

\def\wt{\widetilde}

             %Prime/derivative
      %Double prime/double derivative

%
%  Od Mirana:  Script R s krozcem. Uporabljaj kot $\IR$.
%
\def\IR{{\cal R}\kern-2.7mm\hbox{\raise0.8mm\hbox{\tenrm\char23}}}

\begin	
\cl{\bf THE HOMOTOPY PRINCIPLE IN COMPLEX ANALYSIS:}
\cl{\bf A SURVEY}
\bigskip
\cl{ Franc Forstneri\v c}

\vskip 10mm
\cl{\bf Contents}
\bigskip
\ll{Introduction                                                   \dotfill  $\quad 1$}
\smallskip
\ll{1. The homotopy principle and the Oka principle                \dotfill  $\quad 5$}
\smallskip
\ll{2. The Oka principle for mappings: first examples              \dotfill  $\quad 8$}
\smallskip
\ll{3. Mappings of Stein manifolds into subelliptic manifolds      \dotfill  $\,\,\, 11$}
\smallskip
\ll{4. Removing intersections with complex subvarieties.           \dotfill  $\,\,\, 16$}
\smallskip
\ll{5. Embeddings and immersions of Stein manifolds                \dotfill  $\,\,\, 19$}
\smallskip
\ll{6. Embeddings of open Riemann surfaces in the affine plane     \dotfill  $\,\,\, 22$}
\smallskip
\ll{7. Noncritical holomorphic functions and submersions           \dotfill  $\,\,\, 24$}
\smallskip
\ll{References                                                     \dotfill  $\,\,\, 27$}
\vskip 10mm

%
%
%  IMMERSIONS OF SMOOTH MANIFOLDS
%
%
\centerline{\bf Introduction} \medskip

We say that the {\it homotopy principle} holds for a certain analytic or geometric 
problem if a solution exists provided there are no topological 
(or homotopical, cohomological,$\dots$) obstructions. One of
the principal examples is the theory of smooth immersions developed during
1958--61 by S.\ Smale [Sm1, Sm2] and M.\ Hirsch [Hi1, Hi2]: {\it Immersions of a smooth 
manifold $X$ to affine spaces $\R^N$ of dimension $N>\dim X$ 
can be classified up to regular homotopy 
by their tangent maps}, and hence by vector bundle injections from the 
tangent bundle $TX$ to the trivial bundle $X\times \R^N$. In particular,
an immersion $X\to\R^N$ exists \iff\ $TX$ embeds in $X\times\R^N$. 
If $X$ is an open manifold then the same holds also for $N=\dim X$.

Slightly earlier J.\ Nash [N1, N2] proved that {\it every Riemannian manifold admits 
an isometric immersion into a Euclidean spaces with the flat metric}. In the process
he discovered an important new method for inverting certain classes of non-linear partial 
differential operators by using a suitably modified Newton's iteration to pass from 
solutions of the linearized problem to a solution of the non-linear problem. 
The {\it Nash-Moser-Kolmogorov implicit function theorem} became one of the 
key methods for proving the homotopy principle in problems involving 
underdetermined systems of partial differential equations.

The homotopy principle was investigated even earlier in complex analysis 
where the customary expression for this phenomenon is {\it the Oka principle}.
In 1939 Kiyoshi Oka [Oka] studied the {\it second Cousin problem}: 
Given an open covering $\cU=\{U_j\}$ of a complex manifold $X$ and a collection 
of nowhere vanishing holomorphic functions $f_{ij}\in \cO^*(U_{ij})$ 
satisfying the 1-cocycle condition ($f_{ii}=1$, $f_{ij}f_{ji}=1$,
$f_{ij}f_{jk}f_{ki}=1$), the problem is to find a 
collection of nonvanishing \holo\ functions $f_i\in \cO^*(U_i)$ such that 
$f_i=f_{ij}f_j$ on $U_{ij}=U_i\cap U_j$. (Here $\cO^*$ denotes the 
multiplicative sheaf of nonvanishing holomorphic functions on $X$.)
Oka proved that, if $X$ is a {\it domain of holomorphy} in $\C^n$,
{\it a second Cousin problem can be solved by holomorphic functions 
$f_i$ provided that it can be solved by continuous functions}. 
A modern formulation of Oka's theorem 
is that the Picard group ${\rm Pic}(X)=H^1(X;\cO^*)$ of equivalence 
classes of holomorphic line bundles is isomorphic to the group
$H^2(X;Z)$; the isomorphism is obtained from the cohomology sequence
of the `exponential sheaf sequence' $0\to Z \hra \cO  \to \cO^*\to 1$, 
where the map $\cO\to \cO^*$ is $f\to {\exp(2\pi i f)}$. 

\pn

In 1951 K.\ Stein [Stn] introduced an important class of complex manifolds,
now called {\it Stein manifolds}, on which the algebra of global holomorphic 
functions has similar properties as on domains of holomorphy. It soon 
became clear throught the work of Remmert [Rem] that Stein manifolds can 
be characterized as being biholomorphic to closed complex submanifolds of 
the affine complex spaces. (For more precise results see [Na1, Na2] and [Bis].
For the general theory of Stein manifolds we refer to the monographs [GR],
[GRe] and [H\"or].) H.\ Cartan proved that every coherent analytic sheaf 
on a Stein manifold is generated by global sections and has vanishing cohomology 
groups in all dimensions $q\ge 1$ (Theorems A and B). Thus every analytic
problem on a Stein manifold whose only obstruction lies in such a group is solvable. 

An equivalent formulation of Oka's result on the second Cousin problem is that
two holomorphic line bundles on a Stein manifold are holomorphically isomorphic 
provided they are isomorphic as topological complex vector bundles. 
This problem has an immediate extension to vector bundles of rank $q>1$.
The holomorphic equivalence classes of such bundles are represented by the 
cohomology group $H^1(X;\cG_q)$ with coefficients in the non-abelian sheaf 
of (germs of) holomorphic maps $X\to GL_q(\C)$. Cartan's theory does not apply 
directly to such sheaves and one must in some sense linearize the problem.  
This was done by H.\ Grauert in seminal papers [Gra1, Gra2] from 1957-58
in which he proved that {\it for complex vector bundles over Stein manifolds, 
holomorphic classification agrees with topological classification}.
Grauert reduced the problem to the homotopy classification of sections of 
holomorphic fiber bundles with homogeneous fibers, showing that the inclusion of 
the space of \holo\ sections into the space of continuous sections 
is a {\it weak homotopy equivalence}. For expositions of Grauert's theory
see [Ca], [Ram] and [Lei]. An equivariant version of Grauert's theorem 
was proved in [HK]. A different proof and extension to $1$-convex manifolds
was given in [HL1, HL2]. Recently the result was extended 
to $1$-convex complex spaces [LV].

Progress during the 1960's brought improvements and extensions of the 
Hirsch--Smale theory in real geometry and of Grauert's theory in complex geometry. 
Phillips [Ph1, Ph2, Ph3] showed that the homotopy principle, analogous to the Hirsch-Smale 
theory of immersions, holds for smooth submersions and foliations of open manifolds.
Forster [Fs1, Fs2] applied the Oka-Grauert principle to study holomorphic embeddings 
of Stein manifolds in low dimensional affine spaces. Forster and Ramspott [FRa] proved 
the Oka principle in the problem of holomorphic complete intersections. 
In another direction, Gunning and Narasimhan [GN] (1967) constructed noncritical 
holomorphic functions on any open Riemann surfaces.

The homotopy principle in real differential topology in geometry 
was revolutionized by Mikhael Gromov in the period 1967--73. In his 
seminal paper [Gr1] Gromov presented the method of {\it convex integration 
of differential relations} which unified many seemingly unrelated 
geometric results (in particular the Smale-Hirsch theory of immersions
and Phillips's result on submersions).  Gromov's methods 
initiated rapid progress, and new examples which fit into this
framework are being found even today. For a comprehensive survey of this
and other methods to prove the homotopy principle we refer to Gromov's monograph 
[Gro2, 1986]. During recent years other monographs appeared on this topic, 
see e.g.\ [Sp] and [EM]. The convex integration method, together 
with other methods for solving global problems such as the {\it removal 
of singularities} (Gromov and Eliashberg [GE, Gro2]), {\it continuous sheaves} 
[Gro2], {\it inversions of differential operators} (Nash [N1, N2], 
Hamilton [Ham], Gromov [Gro2]), today provides one of the cornestones of 
differential topology and geometry.

In this paper we survey the homotopy principle in complex analysis and geometry,
drawing parallels with the real geometry where appropriate. Results of this type are 
commonly referred to (as instances of) the {\it Oka principle} when the underlying 
manifold is Stein. To our knowledge this notion has never
been precisely defined, or at least there is no universal agreement on 
what the definition should be. In the monograph [GRe] of Grauert and Remmert 
one finds on p.\ 145 the following formulation: 
{\it Analytic problems (on Stein manifolds) which can be cohomologically 
formulated have only topological obstructions}. If `cohomologically' is interpreted
in the sense that the obstruction to a given problem lies in a cohomology group with 
coefficients in a coherent analytic sheaf then this is just Cartan's Theorem B. 
The Oka-Grauert theory goes a step further by reducing a holomorphic problem to a 
problem in homotopy theory. Hence one is tempted to include in the Oka principle 
all those analytic problems on Stein manifolds which can be homotopically formulated. 
There is a serious limitation to such attempts since certain analytic problems have 
no solution due to {\it hyperbolicity} (Picard's theorems, Kobayashi hyperbolicity, etc.). 

What is then a sensible notion of the Oka principle which would adequately 
cover the known results? It's probably impossible to find one. For the purposes
of this paper we adopt the convention that
\medskip
\centerline{ \it The Oka principle $=$ the homotopy principle in complex analysis.}
\smallskip
\ni We give precise definitions, conforming to Gromov's [Gro2], in Section 1.

It is not surprising that some of the most powerful  methods to prove the homotopy 
principle in the smooth category do not extend to the holomorphic category. The absence of 
partitions of unity can be substituted to a large extent by Cartan's theory and the 
$\dibar$-methods. A  more serious problem is that boundary 
values completely determine holomorphic objects. This disqualifies the convex integration 
method which is based on extending a solution by induction over the skeleta of a CW-complex. 
Fortunately some of the other methods mentioned above remain 
applicable in Stein geometry. This holds for the {\it elimination of singularities}  
which reduces differential conditions to essentially algebraic conditions along 
submanifolds (or complex subvarieties), as well as for the Nash-Moser theory 
of inverting certain differential operators. 

In 1970's it became clear that progress on many questions in Stein geometry
essentially depended on extending the Oka-Grauert principle to sections of 
more general types of holomorphic fiber bundles, and even of non-locally trivial submersions. 
A crucial contribution was made by Henkin and Leiterer [HL1, HL2] 
who reproved Grauert's theorem using the `bumping method', thereby 
localizing the approximation problems which arise
in the construction of global sections. This turned out to be a key 
point which opened the way to extensions. The potential was realized by 
Gromov in 1989 [Gro3] who introduced the concept of a 
{\it dominating spray} as a replacement of the exponential map 
in the linearization and patching problems which appear in the Oka-Grauert theory.
The presence of a dominating spray on the fiber of a \holo\ bundle
over a Stein manifold implies the Oka principle for its sections.
Results on this topic are given in Section 3.

In Section 4 we look at the question of {\it removing intersections} of 
holomorphic maps from Stein source manifolds with closed complex subvarieties 
of the (not necessarily Stein) target manifolds. Progress in 
this direction was made possible by the techniques developed to 
prove Gromov's Oka principle mentioned above. In the classical case 
of {\it complete intersections} the Oka principle was proved in 1967 by 
Forster and Ramspott [FRa].

In Section 5 we survey the results on {\it embeddings and immersions
of Stein manifolds into affine spaces}. In 1956 R.\ Remmert proved [Rem]
that every Stein manifold $X^n$ admits a proper \holo\ embedding in $\C^{2n+1}$
and immersion in $\C^{2n}$ (see also Narasimhan [Na1, Na2] and Bishop [Bis]). 
For embeddings of smooth manifolds $X^n\to\R^N$ the general minimal dimension 
is $N=2n$. However, an $n$-dimensional Stein manifold $X^n$ behaves
in the sense of homotopy theory as being at most real $n$-dimensional. After the 
initial work of Forster [Fs1, Fs2] the optimal embedding dimension $N=[3n/2]+1$ 
was conjectured by Gromov and Eliashberg in 1971 [GE] and proved in 1992 [EG] 
(for odd $n$ the proof was completed in [Sch]). The problem of embedding 
open Riemann surfaces to $\C^2$ is still open and we survey it in Section 6. 

In Section 7 we mention some of the results from the recent work [F8] on the 
existence of holomorphic submersions of Stein manifolds to complex Euclidean spaces. 

The  survey is not comprehensive in any way and the choice of material reflects 
personal tastes and limitations of its author. Among the many topics which are not 
discussed we mention the existence and homotopy classification of Stein structures 
on even dimensional smooth manifolds; see Eliashberg's paper [El] and the 
monographs by Gromov [Gro2] and Gompf and Stipsicz [GSt]. 
I apologize to all authors whose contributions to this field may have
been unjustly left out.

\beginsection 1. The homotopy principle and the Oka principle.

We begin by recalling from [Gro2] the notion of a {\it differential relation}. 
Consider a smooth submersion $h\colon Z\to X$ between smooth real manifolds.
Let $Z^{(r)}$ denote the space of $r$-jets of (germs of) smooth sections 
$f\colon X\to Z$ for $r=0,1,2,\ldots$. The $0$-jet of $f$ at $x\in X$ 
is its value $f(x)\in Z_x=h^{-1}(x)$. The $r$-jet $j^r_x(f) \in Z^{(r)}$
is determined in local coordinates near $x\in X$ resp.\ $f(x)\in Z$ by the 
partial derivatives of $f$ of order $\le r$ at $x$. 

We have natural projections $p^r\colon Z^{(r)}\to Z$ and 
$p^s_r \colon Z^{(s)}\to Z^{(r)}$ for $s>r\ge 0$, where $Z^{(0)}=Z$.
The jet bundles $Z^{(r)}$ carry natural smooth structures, as well 
as affine structures in fibers (see [Gro2] for more details). 
When $X$ and $Z$ are complex manifolds and $h\colon Z\to X$ is a
holomorphic submersion, we shall denote by $Z^{(r)}$ the space of 
$r$-jets of holomorphic sections $f\colon X\to Z$. 

Note that for every section $g\colon X\to Z^{(r)}$ we get a  
corresponding `base point' section $f=p^r(g)\colon X\to Z$.
In general $g$ need not equal $j^r(f)$;  when $g=j^r(f)$ we say that 
the section $g$ is {\it holonomic}.

%
% DIFFERENTIAL RELATION
%
\proclaim 1.1 Definition: {\rm [Gro2, p.\ 2])}
A {\bf differential relation} of order $r$ is a subset 
$\cR \subset Z^{(r)}$ of the $r$-jet bundle $Z^{(r)}$. 
A $\cC^r$ section $f\colon X\to Z$ is said to {\it satisfy}
(or to be a solution of) $\cR$ if $j^r(f) \colon X\to Z^{(r)}$
has values in $\cR$ (i.e., $j^r_x(f)$ belongs to the fiber
$\cR_{f(x)}=(p^r)^{-1}(f(x))$ of $\cR$ over the point $f(x)\in Z$).

The relation $\cR\subset Z^{(r)}$ is said to be {\it open} (resp.\ {\it closed})
when $\cR$ is an open (resp.\ closed) subset of the jet bundle $Z^{(r)}$.
Natural examples of closed relations which arise in geometric problems 
are unions of submanifolds (or subvarieties) of the jet bundle $Z^{(r)}$,
and open relations as complements of submanifolds (or subvarieties). 
Differential equations are examples of closed differential relations.

%
%  H-PRINCIPLE
%
\smallskip \ni \bf 1.2 Definition.  \sl
(a) {\rm ([Gro2, p.\ 3])} 
Let $r$ be a nonnegative integer and let $s\in \{r,r+1,\ldots,\infty\}$.
We say that solutions of class $\cC^s$ of a differential relation 
$\cR\subset Z^{(r)}$ satisfy the {\bf basic h-principle} if every 
continuous section $\phi_0\colon X\to \cR$ is homotopic through sections
$\phi_t \colon X\to \cR$ $(t\in [0,1])$ to a holonomic section
$\phi_1=j^r(f)$ for some $\cC^s$ section $f\colon X\to Z$.  

\ni (b) Assume that $h\colon Z\to X$ is a \holo\ submersion. 
We say that sections $X\to Z$ of $h$ satisfy the {\bf basic Oka principle} 
if every continuous section is homotopic to a \holo\ section. 
(For the parametric Oka principle see Definition 2.1.)

\ni (c) {\rm ([Gro2, p.\ 66]; assumptions as in (b).)} 
For $r\ge 1$ we say that a differential relation $\cR \subset Z^{(r)}$ 
satisfies the {\bf holomorphic h-principle} if every {\rm \holo} section 
$\phi_0\colon X\to \cR$ is homotopic through \holo\ sections 
$\phi_t\colon X\to \cR$ to a holonomic \holo\ section 
$\phi_1=j^r(f) \colon X\to\cR$ (where $f\colon X\to Z$ is a \holo\ section
of $h\colon Z\to X$). We say that $\cR$ satisfies the {\bf basic Oka principle} 
if every {\rm continuous} section $\phi_0\colon X\to \cR$ is homotopic through 
a family of continuous sections of $\cR$ to a \holo\ holonomic section of $\cR$.

%
%  REMARKS
%
\smallskip\ni\bf  1.3 Remarks. \rm
(a) In the $\cC^s$-smooth case one usually takes the {\it fine\/} $\cC^s$ topology on 
the space of $\cC^s$ sections $X\to Z$. In the holomorphic case one must use 
the weaker compact-open topology to obtain meaningful results.

\smallskip \ni
(b) One can introduce more refined notions
such as the {\it parametric h-principle}, the {\it h-principle with approximation}
(or {\it interpolation}), the {\it relative h-principle}, etc. We refer the reader 
to [Gro2]. In section two we introduce some of these notions in 
the holomorphic case (for relations of order zero).

\smallskip \ni
(c) The problem of deforming a continuous section
$\phi_0\colon X\to \cR$ to a \holo\ holonomic section 
can be treated in two steps: 
\item{-} first deform $\phi_0$ through continuous sections of $\cR$ to a
\holo\ section $\phi_1\colon X\to \cR$ (the ordinary Oka principle
for sections $X\to \cR$);
\item{-} deform a (non-holonomic) \holo\ section $\phi_1\colon X\to \cR$
through a homotopy of \holo\ sections $\phi_t\colon X\to \cR$ 
$(t\in [1,2])$ to a \holo\ holonomic section $\phi_2=j^r(f)$
(the holomorphic h-principle).

%
%  EXAMPLES
%
\smallskip\ni\bf 1.4 Examples.  (a)  Mappings $X\to Y$. \rm
An open differential relation of order zero is specified by an open subset 
$\Omega\subset Y$, and the h-principle requires that every continuous map 
$X\to\Omega$ is homotopic to a smooth (real-analytic, holomorphic) map through a 
homotopy with range in $\Omega$. For smooth maps this follows 
from Whitney's approximation theorem. The problem is highly nontrivial in 
the holomorphic case (Sections 2 and 3).

\smallskip \ni \bf (b) Smooth immersions. \rm 
Let $X$ be a smooth manifold. A  map 
$f=(f_1,\ldots,f_q)\colon X\to\R^q$ is an immersion 
if its differential $df_x\colon T_xX \to T_{f(x)}\R^q\simeq \R^q$ is a {\it injective} 
for every $x\in X$. The pertinent differential relation (of order one) consists of 
all points $(x,y,\l)$ where $x\in X$, $y\in \R^q$ and $\l\in \Hom(T_x X,\R^q)$ 
with $\l$ {\it injective}. Clearly the value $f(x)$ is unimportant due to translation 
invariance, and we can reduce the problem to the relation $\cR$ whose sections are 
{\it injective vector bundle maps} $TX\to X\times \R^q$ from the tangent bundle of 
$X$ into the trivial bundle $X\times\R^q$. (Alternatively, we can consider the relation
whose sections are $q$-tuples of differential 1-forms $\theta=(\theta_1,\ldots,\theta_q)$ 
on $X$ which together span the cotangent space $T_x^* X$ at each point $x\in X$.)  
The h-principle of Smale [Sm1, Sm2] and Hirsch [Hi1, Hi2] asserts 
that if either $q>\dim X$, or if $q=\dim X$ and $X$ is open, then  
{\it the regular homotopy classes of smooth immersions $X\to\R^q$ are in 
one-to-one correspondence with the homotopy classes of vector bundle injections 
$TX\to X\times\R^q$}. In particular, an immersion $X\to\R^q$ exists
\iff\ the cotangent bundle $T^*X$ is generated by $q$ sections. 

\smallskip \ni\bf (c) Holomorphic immersions. \rm
The Oka principle for holomorphic immersions $X\to\C^q$ of Stein manifolds 
to affine spaces of dimension $q>\dim X$ was proved by Eliashberg and Gromov [Gro2] 
(Section 5 below). The problem is open in the critical dimension $q=\dim X$ 
except for a positive result in dimension $n=1$ due to Gunning and Narasimhan [GN]. 
The Oka principle also holds for {\it relative immersions} (maps 
$g\colon X\to\C^n$ such that $f=b\oplus g\colon X\to\C^{m+n}$ is a \holo\ 
immersion, where $b\colon X\to \C^m$ is a fixed holomorphic map). 
 
\smallskip
\ni \bf (d) Smooth submersions. \rm
These are smooth  maps $X\to Y$ of rank equal to $\dim Y$ at each 
point of $X$ (hence $\dim Y\le \dim X$). The tangent map of a submersion 
$X^n\to  \R^q$ induces a surjective vector bundle map $TX\to X\times\R^q$.
The homotopy principle  due to Phillips [Ph1, Ph3] asserts that 
{\it for any open manifold $X$, the regular homotopy classes of submersions 
$X\to\R^q$ are in one-to-one correspondence with surjective vector bundle maps 
$TX\to X\times \R^q$.} In particular, a submersion $X\to\R^q$
exists if the tangent bundle $TX$ admits a trivial subbundle
of rank $q$. (See also Gromov [Gro2, pp.\ 26, 53].) 
In fact the h-principle holds for smooth submersions of open manifolds
to arbitrary manifolds. It also holds for smooth maps $X\to Y$ 
of constant rank $k$ [Ph2] or rank $\ge k$ (Feit [Fe], [Gro2, p.\ 27]). 

\smallskip
\ni \bf (e) Holomorphic submersions. \rm
In 1967 Gunning and Narasimhan proved that {\it every open Riemann surface admits a 
holomorphic function without critical points} [GN]. Very recently it was proved 
in [F8] that the same holds on every Stein manifold $X$. Moreover, 
{\it \holo\ submersions $X\to \C^q$ satisfy the Oka principle when $q<\dim X$} 
(Section 7 below). In the maximal dimension $q=\dim X$ the problem is still open: 
{\it Does every Stein manifold $X^n$ with trivial tangent bundle admit a 
locally biholomorphic map $f\colon X\to\C^n$?} 

\smallskip
\ni \bf (f) Foliations. \rm
A foliation $\cF$ of rank $q$ of an $n$-dimensional manifold $X$ is determined by 
an integrable rank $q$ subbundle $E \subset TX$ whose fiber $E_x$ is the tangent 
space to the leaf of $\cF$ through $x\in X$. The corresponding homotopy 
principle was proved for smooth open manifolds by Phillips [Ph2],
and for closed manifolds by Thurston [Th1, Th2] 
(see also [Gro2], p.\ 102 and p.\ 106]): 
{\it If a smooth subbundle $E\subset TX$ has a trivial normal bundle 
$TX/E$ then $E$ is homotopic to an integrable smooth subbundle in $TX$}
(which therefore determines a smooth foliation on $X$). The same holds 
if $N=TX/E$ admits locally constant transition functions. 
In [F8] the analogous results are proved for holomorphic foliations on 
Stein manifolds. The h-principle fails 
for real-analytic foliations on {\it closed\/} manifolds since 
{\it no closed simply connected real-analytic manifold admits a 
real-analytic foliation of codimension one} (Haefliger [Hae]). 
For example, the seven-sphere $S^7$ admits a smooth codimension 
one foliation but no real-analytic one.

\smallskip \ni
{\bf (g) Totally real immersions and embeddings.}
Let $S$ be a smooth manifold and $X$ a complex (or almost complex) manifold.
An immersion $f\colon S\to X$ is {\it \tr\/} if for each $p\in S$
the image $\Lambda_p = df_p(T_pS) \subset T_{f(p)}X$ is a \tr\ linear subspaces
of $T_{f(p)} X$, i.e, $\Lambda_p\cap J(\Lambda_p)=\{0\}$ where 
$J \in {\rm End}(TX)$ is the almost complex structure  on $X$.
The pertinent differential relation is the set of tripples
$(p,x,\l)$ where $p\in S$, $x\in X$ and $\l \in \Hom(T_p S, T_x X)$ 
is an injective $\R$-linear map whose image is a \tr\ subspace of $T_x X$.
When $X=\C^n$ we can reduce the problem to a relation whose sections are injective 
$\R$-linear vector bundle maps $\iota\colon TS\to S\times \C^n$  
with $\iota(T_pS) \cap i\iota(T_p S) =\{0\}$ for all $p\in S$. The homotopy principle 
holds for \tr\ immersions and embeddings of any smooth manifold into any complex 
manifold. An important point is that the {\it Whitney trick can be performed
through totally real submanifolds}. See [Gro2, EH, F1, Au], and for dimension 
two also [F2, F7]. 

\smallskip\ni 
{\bf (h) Isometric immersions and embeddings.}
The fundamental work of Nash [Na1, Na2] initiated a rich and complex theory.
An interesting feature of Nash's discovery is that every smooth Riemannian 
manifold $(X^n,g)$ admits local isometric immersions of class $\cC^1$ already 
in $\R^{n+1}$ (codimension one!), but $\cC^\infty$ isometric immersions 
require larger codimension. For $q\ge {1\over 2}(n+2)(n+3)$ 
{\it free isometric immersions $X^n\to\R^q$ satisfy the homotopy principle} 
[Gro2, p.\ 12]. Important results in this field were obtained 
by R.\ Greene and collaborators [Gn, GW, GSh]. We refer to 
[Gro2] for further results and references.

%
%
%  SECTION 2
%
%
\beginsection 2. The Oka principle for mappings: first examples.

Let $X$ and $Y$ be complex manifolds. We say that mappings 
$X\to Y$ satisfy the {\bf basic Oka principle} if every continuous map
is homotopic to a \holo\ map. For mappings between Riemann surfaces
a complete answer on the validity of the Oka principle was given
by Winkelmann [Wi]. For manifolds of higher dimension only partial 
results are known.

A compact set $K$ in a Stein manifold $X$ is {\it holomorphically convex}
if for every point $p\in X\bs K$ there is an $f\in \cO(X)$ such that
$f(p_0) > \sup_{x\in K} |f(x)|$.
Let $P$ be a compact Hausdorff space (the parameter space) and $P_0\subset P$ 
be a compact subset which is a strong deformation retraction of some \nbd\ of 
$P_0$ in $P$. We shall consider the following stronger versions of the 
Oka principle (see the notion of $Ell_\infty$ fibrations in [Gro3], as well as 
Theorem 1.5 in [FP2].)

\smallskip\ni\bf 2.1 Definition. \sl
(a) Maps $X\to Y$ satisy the {\bf parametric Oka principle} if 
for any continuous map $f\colon X\times P\to Y$ such that 
$f_p=f(\cdotp,p)\colon X\to Y$ is \holo\ for each $p\in P_0$ 
there is a homotopy of continuous maps $f^t\colon X\times P\to Y$ 
$(t\in [0,1])$ such that (i) $f^0=f$, (ii) the map 
$f^1_p :=f^1(\cdotp,p) \colon X\to Y$ is \holo\ for each $p\in P$, and 
(iii) $f^t_p=f_p$ for each $p\in P_0$ and $t\in [0,1]$.

\ni (b) Let $X$ be a Stein manifold and let $\rho$ be a metric on 
the complex manifold $Y$ inducing the manifold topology. 
Maps $X\to Y$ satisfy the {\bf parametric Oka principle with approximation} 
if for every compact holomorphically convex subset $K\subset X$ and for every 
$f\colon X\times P\to Y$ as in (a) such that all maps $f_p$ $(p\in P)$ are 
\holo\ in an open \nbd\ of $K$ in $X$, there is for every $\e>0$ a homotopy 
$f^t \colon X\times P\to Y$ satisfying (a) and also $\rho(f^t(x,p), f(x,p)) <\e$ 
$(x\in K,\ p\in P)$.

\ni (c) Maps $X\to Y$ satisfy the {\bf parametric Oka principle with interpolation} 
if for every closed complex subvariety $X_0\subset X$ and for every 
$f\colon X\times P\to Y$ as in (a) such that all maps $f_p$ $(p\in P)$ are \holo\ in 
an open \nbd\ of $X_0$ in $X$, there is for every $r\in \N$ a homotopy 
$f^t\colon X\times P\to Y$ as in (a) such that $f^t_p$ 
agrees with $f_p$ to order $r$ on $X_0$ for every $t\in [0,1]$ and $p\in P$.
\rm

\smallskip
The above notions clearly extend to sections $f\colon X\to Z$ of a \holo\ 
submersion $h\colon Z\to X$ (see [Gro3] or [FP2]).

\smallskip\ni\bf 2.2 Remark. \rm
The parametric Oka principle for maps $X\to Y$ implies that the inclusion 
$$ 
	\iota \colon {\rm Holo}(X;Y) \hra {\rm Cont}(X;Y)  \eqno(2.1)
$$
of the space of \holo\ maps into the space of continuous maps is a 
{\it weak homotopy equivalence}, i.e., it induces isomorphisms 
of the corresponding homotopy groups of the two spaces which are 
equipped with the compact-open topology [FP1]. (In some papers this is
the definition of the parametric Oka principle.) 
In particular, each connected component of ${\rm Cont}(X;Y)$
contains precisely one connected component of ${\rm Holo}(X;Y)$
which means that (a) every continuous map is homotopic to a \holo\ map, and
(b) every homotopy between a pair of holomorphic maps can be 
continuously deformed to a homotopy  consisting of holomorphic maps.
\endpr

The basic Oka principle can hold for a trivial reason that either of the 
manifolds $X$ or $Y$ is contractible (hence every map is homotopic to constant). 
However, topological contractibility does not necessarily imply the 
parametric (or any other) Oka principle. In all nontrivial situations where 
the Oka principle has been established for all Stein source manifolds $X$, 
it was actually proved in the strongest form (parametric, with interpolation and 
approximation). This is no coincidence since at least the approximation is built 
into all the known proofs.  

We collect positive results on the Oka principle for maps (and sections)
in Section 3 below. In the remainder of this section we give examples 
illustrating the failure of the Oka principle for maps from Stein
source manifolds. In most examples the reason for the failure is 
{\it hyperbolicity\/} of the target manifold. 

\smallskip\ni\bf 2.3 Example. \rm
If either $X$ or $Y$ is contractible then the Oka principle 
for maps $X\to Y$ trivially holds. However, the Oka principle with approximation 
fails already for self-maps of the unit disc $U=\{z\in \C\colon |z|<1\}$: 
If $a\in U\bs \{0\}$ and $0<r<1-|a|$, the translation $f(z)=z+a$ maps 
$\{|z|< r\}$ holomorphically into $U$ and it extends to a smooth a map $U\to U$, 
but it cannot be approximated uniformly on any \nbd\ of the origin $0\in U$ by a 
\holo\ map $g\colon U\to U$ since this would give $g'(0)\approx f'(0)=1$ in 
contradiction to the Schwarz lemma. 

\smallskip \ni\bf 2.4 Example. \rm
This and the next example can be found in [Gro3] (see also [Wi]).
Let $X=\{1<|z|<r\}$ and $Y=\{1<|z|<R\}$ be annuli in $\C$. The space of homotopy 
classes of maps $X\to Y$ equals $\pi_1(Y)=Z$. However, if $1< R<r$ then every 
holomorphic map $X\to Y$ is homotopic to constant and hence the Oka principle fails. 
Furthermore, for any choice of values of $r, R>1$ only finitely many homotopy 
classes of maps $X\to Y$ are represented by holomorphic maps. To see this, observe that 
the infimum of the Kobayashi length of closed curves in $X$ or $Y$ which 
generate the respective fundamental group is positive, and holomorphic maps 
do not increase the length. On the other hand, the Oka principle 
holds for maps $X\to \C\bs \{0\}$ from any Stein manifold $X$ 
(Section 3).

\smallskip \ni\bf 2.5 Example. \rm
The argument in Example 2.4 extends to any Kobayashi hyperbolic 
target manifold. Recall that a complex manifold $Y$ is {\it Kobayashi hyperbolic} if for any 
point $y\in Y$ and tangent vector $v\in T_y Y\in \{0\}$ the set of all numbers 
$\l\in\C$ of the form $f'(0)=\l v$ for some \holo\ map $f\colon U\to Y$, $f(0)=y$, 
is bounded:  $|\l| \le M$ for some $M=M(y,v) < +\infty$. For instance, the twice 
punctured plane $\C\bs \{0,1\}$ is hyperbolic by Picard's theorem. 
If we take as before $X$ to be an annulus then even the basic Oka principle 
fails for maps $X\to Y= \C\bs \{0,1\}$ which is seen by the following 
argument from [Gro3, p.\ 853]. Take a circle $S = \{|z|=\rho\} \subset X$ and 
wrap is sufficiently many times around each of the points $0,1$  by a smooth map 
$f\colon S \to Y$. Since the minimal Kobayashi length of closed curves 
representing a given homotopy class in $\pi_1(Y)$ increases to $+\infty$ 
when we increase the number of rotations around the two punctures, 
the length of $f(S)$ in $Y$ will exceed the length of $S$  in $X$
for any $f$ representing a suitably chosen class in $\pi_1(Y)$. 
Since \holo\ maps do not increase the Kobayashi length, it follows 
that such $f$ is not homotopic to any holomorphic map $X\to Y$. 
In fact only finitely many classes in $\pi_1(Y)$ can be represented 
by holomorphic maps $X\to Y$.

\smallskip \ni\bf 2.6 Example. \rm
The following example, due to J.-P.\ Rosay (private communication), 
is an improvement of Proposition 2.2 in [\v CF]. It shows that holomorphic 
graphs over the unit disc cannot avoid even fairly simple complex 
curves in $\C^2$. Let $(z,w)$ be complex coordinates on $\C^2$.
For $k\in \C$ let 
$$
   \Sigma_k =\{w=0\}\cup \{w=1\} \cup\{w=kz\} \cup \{zw=1\} \subset \C^2.
$$

\proclaim Proposition: {\rm (J.-P.\ Rosay)} There is a 
$k>0$ such that the graph of any \holo\ function 
$f\colon U=\{|z|<1\} \to \C$ intersects $\Sigma_k$. (Indeed this is
true for every sufficiently large $|k|$.) On the other hand, 
for any $k$ there exists a smooth function $U\to \C$ whose 
graph avoids $\Sigma_k$.

This should be compared with Examples 3.4 and 3.5 below on avoiding subvarieties 
of codimension at least two. This is also in strong contrast to the situation 
for {\it holomorphic motions}, i.e., {\it disjoint unions of holomorphic graphs} 
over the disc, which can always be extended to maximal motions according 
to Slodkowski [Sl1, Sl2].

\demo Proof: The last statement is a simple topological exercise.
Suppose now that $f\colon U\to \C\bs\{0,1\}$ is a holomorphic function
omitting $0$ and $1$. Denote by $l$ the length of the circle 
$C=\{|z|=1/2\}$ with respect to the Kobayashi ($=$Poincar\'e)
metric on $U$. Denote by $d$ the Kobayashi distance function
on $\C\bs \{0,1\}$. Let
$$ 
	k_0= \sup\{|\z|\in \C\bs \{0,1\}\colon 
        \inf_\theta d(\zeta,2e^{i\theta}) \le l\}. 
$$
We have $k_0<+\infty$ since $\C\bs\{0,1\}$ is complete hyperbolic.
Observe that the Kobayashi length of $f(C) \subset \C\bs\{0,1\}$
is at most $l$. We consider two cases.

\demo Case 1: There exists a $\theta\in \R$ such that 
$|f(e^{i\theta}/2)|\le 2$. Then for all $\gamma\in \R$
we have $|f(e^{i\gamma}/2)|\le k_0$ by the choice of $k_0$.
Rouch\'e's theorem shows that for every $k> 2k_0$ 
the equation $kz - f(z)=0$ has a solution with $|z|<1/2$,
and at this point the graph of $f$ intersects $\Sigma_k$.

\demo Case 2: For every $\theta\in \R$ we have
$|f(e^{i\theta}/2)| > 2$. Since $f$ has values in $\C\bs\{0,1\}$,
so does $g=1/f$, and the above gives $|g(e^{i\theta}/2)| < 1/2$
for every $\theta$. Rouch\'e's theorem implies that 
$z-g(z)=z-{1\over f(z)}$ has one zero with $|z|<1/2$, which means 
that $zf(z)=1$ has a solution with $|z|<1/2$. At this point
the graph of $f$ intersects $\Sigma_k$.
This completes the proof.

\smallskip \ni\bf 2.7 Example. \rm This example is taken from [FP3]. 
{\it For every $n\in \N$ there exists a discrete subset $P\subset \C^n$
such that the basic Oka principle fails for maps $X\to \C^n\bs P$
from some Stein manifold $X$.} In fact this holds for any discrete set $P$ 
which is {\it unavoidable} in the sense of Rosay and Rudin [RR], i.e., such that
any entire map $\C^n\to\C^n$ of generically maximal rank intersects 
$P$ infinitely often. Alternatively, any entire map $\C^k\to\C^n\bs P$
has rank $<n$ at each point (here $k$ may be different from $n$). 
The same holds for maps $X\to \C^n\bs P$ for any $X$ covered by an
affine space. 

A simple argument shows that any holomorphic map $X\to \C^n\bs P$ of rank 
$<n$ is homotopic to the constant map in $\C^n\bs P$. However, for certain 
$X$ covered by an affine space there exist homotopically nontrivial smooth maps 
$X\to\C^n\bs P$ and hence the Oka principle fails. To obtain such an example
let $n=2$ and $X=(\C\bs \{0\})^3$ (which is universally covered by $\C^3$). 
There is a smooth contraction of $X$ onto the standard torus 
$T^3\subset \C^3$. Let $f\colon T^3\to \C^2\bs P$ be an embedding of $T^3$ 
onto a small hypersurface torus surrounding a point $p_0\in P$. Composing 
$f$ with the contraction $X \to T^3$ we get a nontrivial smooth map 
$X\to \C^2\bs P$ which is not homotopic to any holomorphic map. 

Note that the infinitesimal Kobayashi pseudometric on $Y=\C^n\bs P$ is totally 
degenerate, but the Kobayashi-Eisenmann volume form on $Y$ is nontrivial 
(when $P$ is unavoidable).
\beginsection 3. Mappings of Stein manifolds into subelliptic manifolds.

In this section we present results on the Oka principle for maps 
$X\to Y$ from Stein source manifolds, as well as for section of submersions 
onto a Stein base. Our main references are the papers by
Grauert [Gra1, Gra2], Cartan [Ca], Gromov [Gro3], and [FP1, FP2, FP3, F5].

By definition every Stein manifold admits plenty of holomorphic maps to complex
affine spaces $\C^q$. The basic idea introduced by Gromov [Gro3] is the following.
Suppose that a complex manifold $Y$ admits sufficiently many 
{\it dominating holomorphic maps\/} $s\colon \C^q\to Y$, where the domination property 
means $s$ is a submersion outside a subvariety of $\C^q$. Then there also exist plenty 
of \holo\ maps $X\to Y$ from any Stein $X$. (In some sense the idea is to factor maps 
$X\to Y$ as $X\to\C^q\to Y$.) What is needed in the proofs is a family of dominating
maps $s_y\colon \C^q\to Y$, depending holomorphically on the point $y=s_y(0) \in Y$. 
This leads to the following concept of a {\it dominating spray} introduced by Gromov 
[Gro3]. The notion of a {\it dominating family of sprays} and of 
{\it subelliptic manifolds} was introduced in [F5].

%
%  SPRAY
%
\proclaim 3.1 Definition: 
A {\bf spray} on a complex manifold 
$Y$ is a holomorphic map $s\colon E\to Y$, defined on the total
space of a holomorphic vector bundle $p\colon E\to Y$, such that
$s(0_y)=y$ for every $y\in Y$. The spray is {\bf dominating at $y$} if 
its differential $ds_{0_y}\colon T_{0_y}E\to T_yY$ maps $E_y$ (which is 
a linear subspace of $T_{0_y}(E)$) onto $T_y Y$; it is {\bf dominating} 
if this holds at every point $y\in Y$. A {\bf dominating family of sprays} 
is a collection of sprays $s_j\colon E_j\to Y$ $(j=1,2,\ldots,k)$ 
such that for every $y\in Y$ we have
$$ 
	(ds_1)_{0_y}(E_{1,y}) + (ds_2)_{0_y}(E_{2,y})\cdots 
                     + (ds_k)_{0_y}(E_{k,y})= T_y Y.   \eqno(3.1)
$$
A manifold $Y$ is called {\bf elliptic} if it admits a dominating spray,
and {\bf subelliptic} if it admits a finite dominating family of sprays.

\proclaim 3.2 Theorem: 
{\bf (The Oka principle for maps to subelliptic manifolds.)}
If $X$ is a Stein manifold and $Y$ is a subelliptic manifold then mappings 
$X\to Y$ satisfy the parametric Oka principle with interpolation and approximation. 
Furthermore, the Oka principle holds (in all forms) for sections $X\to Z$ of any 
holomorphic fiber bundle $h\colon Z\to X$ with subelliptic fiber $Z_x=h^{-1}(x)$.

We emphasize that there is no restriction on the structure group of $Z\to X$
(we may use the entire group of holomorphic automorphisms of the fiber).
Theorem 3.2 includes the results of Grauert [Gra1, Gra2] and Gromov [Gro3, Sec.\ 2].
Observe that (sub)ellipticity of a complex manifold eliminates Kobayashi or 
Eisenman hyperbolicity.

Theorem 3.2 is proved constructively by approximation and gluing of holomorphic 
maps to $Y$ (resp.\ of sections $X\to Z$) defined on holomorphically 
convex subsets of $X$. The main steps have been developed in the papers cited above 
and in [FP1, FP2, FP3]. In [F5] the result was proved in the final form as stated here.  
For an extension to sections of subelliptic submersions see Theorem 3.8 below.

We now give examples of sprays and (sub)elliptic manifolds; thus the Oka
principle holds for maps from any Stein manifold to any manifold on this list.
Most of them can be found in [Gro3]. 

\smallskip\ni \bf 3.3 Example: Complex homogeneous manifolds. \rm 
Let $G$ be a complex Lie group which acts holomorphically and transitively 
on a complex manifold $Y$ by \holo\ automorphisms. Let ${\bf g}=T_e G$ 
denote its Lie algebra and $\exp\colon {\bf g}\to G$ the associated exponential map. 
The map $Y\times {\bf g}\to Y$, $(y,t)\to e^t y$, is a dominating spray 
on $Y$ and hence $Y$ is elliptic. The Oka principle for maps 
$X\to Y$ is due to Grauert [Gra1, Gra2].

\smallskip\ni \bf 3.4 Example: Sprays induced by complete vector fields. \rm
(See [Gro3] and [FP1].) 
Let $V_1,\ldots,V_k$ be holomorphic vector fields on $Y$ which are complete 
in complex time. Denote by $\theta_j\colon Y\times\C \to Y$ the flow of $V_j$. 
The superposition of these flows (in any order) gives a spray 
$s\colon Y\times \C^k\to Y$ which is dominating at a point $y\in Y$ if the 
vectors $V_1(y),\ldots,V_k(y)$ span the tangent space $T_y Y$. For example, if 
$A\subset \C^n$ is an algebraic subvariety of complex codimension at least two 
then the complement $Y=\C^n\bs A$ admits a dominating polynomial spray of this kind.
(This fails for most complex hypersurfaces $A$, see Example 2.6.)
A complex Lie group admits sprays of this kind induced by left (or right) 
invariant vector fields spanning the Lie algebra.

\smallskip\ni \bf 3.5 Example: Complements of projective subvarieties. \rm
{\it If $A$ is a closed complex ($=$algebraic) subvariety
of complex codimension at least two in the complex projective space
$\CP^n$ then the manifold $Y=\CP^n\bs A$ is subelliptic. The same holds 
if we replace $\CP^n$ by a complex Grassmanian.} The proof proceeds as follows 
(see Proposition 1.2 in [F5]; the idea can be found in [Gro3]). 
Removing a complex hyperplane $L$ from $\CP^n$ we are left with
$\C^n\bs A$ which admits a dominating algebraic (polynomial) spray
defined on a trivial bundle $E=(\C^n\bs A)\times\C^k \to \C^n\bs A$ (Example 3.4). 
Let $[L]=\cO_{\CP^n}(1)$ denote the line bundle on $\CP^n$ 
determined by the divisor of $L$ (the `hyperplane section bundle').
For sufficiently large $m>0$ the spray $s$ extends to an algebraic
spray $\wt s\colon E\otimes [L]^{-m} \to \CP^n$ which is dominating 
over $\CP^n\bs (L\cup A)$. Repeating this with $n+1$ hyperplanes in 
general position we obtain a finite dominating family of sprays on $Y$. 
The bundles of these sprays are nontrivial (in fact negative) and hence 
we cannot combine them into a single dominating 
spray as we did with the flows $\theta_j$ in Example 3.4. It is not known 
whether $\CP^n \bs A$ is elliptic for every such $A$. 

\smallskip\ni \bf 3.6 Example: Matrix-valued maps with nonzero determinant. \rm
Let $X$ be a Stein manifold and $g_1,\ldots,g_k\colon X\to \C^n$
holomorphic maps (with $1\le k< n$) such that the vectors 
$g_1(x),\ldots,g_k(x)\in\C^n$ are $\C$-linearly independent 
for every $x\in X$. The problem is to find \holo\ maps 
$g_{k+1},\ldots,g_{n} \colon X\to \C^n$ such that the matrix
$g(x)=(g_1(x),\ldots,g_n(x))$ with columns $g_j(x)$ satisfies  
$\det g(x) \ne 0$ (or even $\det g(x)=1$) for every $x\in X$. 
The Oka principle holds in this problem (in all forms); in particular, 
{\it a \holo\ solution exists provided there exists a continuous solution.}
To see this we consider the manifold 
$$ \eqalign{
	Z &= \{(x,v_{k+1},\ldots,v_n)\colon x\in X,\ v_j\in \C^n\ 
	{\rm for\ }j=k+1, \ldots, n,   \cr
	&\qquad\qquad 
	\det(g_1(x),\ldots,g_k(x),v_{k+1},\ldots,v_n)\ne 0 \} \cr}
$$
with the projection $h\colon Z\to X$ onto the first factor. 
A solution to the problem is a section $X\to Z$ of this fibration.
The Oka principle follows from the observation that 
$Z\to X$ is a \holo\ fiber bundle whose fiber $GL_{n-k}(\C)\times \C^{k(n-k)}$
is a complex Lie group. For maps into $SL_n(\C)$ is suffices to 
divide one of the columns by the (nonvanishing) determinant function.
\smallskip

Theorem 3.2 extends to sections of {\it subelliptic submersions} which 
we now introduce. Let $h\colon Z\to X$ be a \holo\ submersion onto $X$.
For $U\subset X$ we write $Z|_U=h^{-1}(U)$. For $z\in Z$ we denote by $VT_z Z$ 
the kernel of $dh_z$ (which equals the tangent space to the fiber of $Z$ at $z$) 
and call it the {\it vertical tangent space} of $Z$ at $z$.
The space $VT(Z)\to Z$ with fibers $VT_z Z$ is a holomorphic 
vector subbundle of the tangent bundle $TZ$. 

If $p\colon E\to Z$ is a \hvb\ we denote by $0_z\in E$ the 
base point in the fiber $E_z=p^{-1}(z)$. At each point 
$z\in Z$ ($=$the zero section of $E$) 
we have a natural splitting $T_{0_z}E = T_z Z\oplus E_z$.

%
%
%  DEFINITION OF A SPRAY; ELLIPTICITY
%
%
\proclaim 3.7 Definition:  {\rm [Gro3, sec.\ 1.1.B]}
A {\bf spray} associated to a \holo\ submersion $h\colon Z\to X$ 
(an $h$-spray) is a triple $(E,p,s)$, where $p\colon E\to Z$ is 
a \hvb\ and $s\colon E\to Z$ is a \holo\ map such that 
for each $z\in Z$ we have $s(0_z)=z$ and $s(E_z) \subset Z_{h(z)}$.
The spray $s$ is {\rm dominating} at the point $z\in Z$ if the derivative 
$ds \colon T_{0_z} E \to T_z Z$ maps $E_z$ (which is a linear
subspace of $T_{0_z} E$) surjectively onto $VT_z Z=\ker dh_z$. 
The submersion $h\colon Z\to X$ is called {\bf  subelliptic} if each 
point in $X$ has an open \nbd\ $U\subset X$ such that $h\colon Z|_U \to U$ 
admits finitely many $h$-sprays $(E_j,p_j,s_j)$ for $j=1,\ldots,k$ satisfying  
$$ 
	(ds_1)_{0_z}(E_{1,z}) + (ds_2)_{0_z}(E_{2,z})\cdots 
                     + (ds_k)_{0_z}(E_{k,z})= VT_z Z            \eqno(3.2)
$$
for each $z\in Z|_U$. A collection of sprays satisfying (3.1) is said 
to be dominating at $z$. A submersion $h$ is {\bf elliptic} if the above 
holds with $k=1$.

Comparing with Definition 3.1 we see that a spray on a manifold $Y$ is the 
same thing as a spray associated to the trivial submersion $Y\to point$. 
By definition every elliptic submersion is also subelliptic, 
but the converse is not known. A \holo\ fiber bundle $Z\to X$ is (sub)elliptic 
\iff\ the fiber has this property (since a spray on $E$ induces an $h$-spray on 
the product bundle $h\colon U\times E\to U$).

\proclaim 3.8 Theorem: 
If $h\colon Z\to X$ is a subelliptic submersion onto a Stein manifold $X$ 
then sections $f\colon X\to Z$ satisfy the parametric Oka principle with 
interpolation and approximation.

For elliptic submersions Theorem 3.8  coincides with Gromov's Main Theorem 
in [Gro3, Sec.\ 4.5]. The result is proved in [FP1] for fiber bundles with 
elliptic fibers, in [FP2] for elliptic submersions but without interpolation,  
in [FP3] for elliptic submersions with interpolation, and the extension to 
subelliptic submersions is obtained in [F5]. A version of the Oka principle for 
multi-valued sections of ramified holomorphic maps can be found in [F6].
F.\ L\'arusson [L\'ar] explained this result from the homotopy theory point of view.

\smallskip\ni\bf Historical comments. \rm
The so-called {\it classical case} of Theorem 3.8 (for sections of principal fiber 
bundles with complex homogeneous fibers) is due to Grauert theorem  [Gra1, Gra2, Ca, Ram]. 
A very important develoment which opened the way to generalizations was the 
work of Henkin and Leiterer in 1984 [HL1] (published in 1998 [HL2]) 
where they introduced the {\it bumping method\/} to this problem.
Its main advantage over the original method of Grauert is that one only needs the Oka-Weil
approximation theorem for sections of $Z\to X$ over small subsets of the base 
manifold $X$. The second  main contribution was made by Gromov [Gro3] who replaced 
the exponential map on the fiber (which was used in Grauert's proof to linearize 
the problems and to patch local sections) by the more flexible notion of 
a dominating spray. The idea of using several sprays instead of one is  
implicitly present in [Gro3], but the condition which we call 
subellipticity was not formulated there explicitly. 

\smallskip\ni\bf 3.9 Example: Avoiding subvarieties with algebraic fibers. \rm
Let $h\colon Z\to X$ be a fiber bundle with fiber $Z_x=h^{-1}(x)\simeq \CP^n$ 
or a complex Grassmanian. Assume that $A$ is a closed complex subvariety of $Z$ 
whose fiber $A_x=A\cap Z_x$ has complex codimension at least 
two in $Z_x$ for every $x\in X$. ($A_x$ is algebraic by Chow's theorem.)
Then the restricted submersion $h\colon Z\bs A \to X$ is subelliptic and hence 
Theorem 3.8 applies (Proposition 1.2 (b) in [F5]). If $Z$ is 
obtained from a \holo\ vector bundle $E\to X$ by taking projective closure of
each fiber (i.e., $Z_x \simeq \CP^n$ is obtained by adding to $E_x\simeq \C^n$ 
the hyperplane at infinity $\Lambda_x\simeq \CP^{n-1}$ for every $x\in X$), 
the restricted submersion $h\colon E\bs A\to X$ is even elliptic 
(Corollary 1.8 in [FP2]). The Oka principle for such submersions is 
used in the constructions of proper holomorphic immersion and embeddings of 
Stein manifolds in affine spaces of minimal dimension (Section 5).

\smallskip\ni\bf 3.10 Ellipticity versus subellipticity. \rm
By definition every elliptic complex manifold is also subelliptic.
It is not known whether there exist subelliptic manifolds which are not
elliptic. Natural candidates for a possible counterexample are the complements 
$Y=\CP^n\bs A$ of generic algebraic subvarieties $A\subset\CP^n$ of codimension 
at least two (and of sufficiently large degree); see Example 3.5. 
Such $Y$ admits a finite dominating family of algebraic sprays defined on 
{\it negative} line bundles over $\CP^n$, but we don't see how to 
obtain a dominating spray. (See [F5] for more.)

\smallskip\ni\bf 3.11 Problem. \rm   It is not known whether 
the Oka principle holds for maps from Stein manifolds into the
complement $Y=\C^n\bs K$ of any infinite compact set $K\subset \C^n$
(this is unknown even if $K$ is a closed ball). Such complements have no 
Kobayashi-Eisenman hyperbolicity. If $K$ is convex then every point 
$y\in Y=\C^n\bs K$ is contained in a Fatou-Bieberbach domain 
$\Omega_y\subset Y$. However, it is not known whether $Y=\C^n\bs K$ 
admits any nontrivial sprays. (It is easily seen that there are no sprays 
$s\colon Y\times \C^N\to Y$ from a trivial bundle.) A good test case 
for the validity of the Oka principle might be the \holo\ map
$$ 
	{\rm SL}(2,\C)\to\C^2\bs\{0\}, \qquad
        \pmatrix{\a&\b\cr \g&\d} \to (\a,\b) \qquad  (\a\d-\b\g=1).
$$
This map is clearly homotopic to a smooth map into $Y=\{z\in \C^2\colon |z|>1\}$
but it is unknown whether it is homotopic to a \holo\ map to $Y$. 
(This problem has been mentioned in [FP3].)

The main problem in this connection is the following. Suppose that 
$C\subset B \subset \C^m$ is a pair of compact convex sets. 
Let $K$ be a closed ball in $\C^n$ and let $f\colon \wt C\to \C^n\bs K$ 
be a \holo\ map on a \nbd\ of $C$ whose range avoids $K$. Is it possible to 
approximate $f$ uniformly on $C$ by a \holo\ map $\wt f \colon \wt B \to \C^n\bs K$
defined on a \nbd\ of $B$~?

\smallskip\ni\bf 3.12 Problem. \rm 
For any Stein manifold $Y$ Theorem 3.2 has the following converse:
{\it If the Oka principle holds for maps $X\to Y$ from any Stein manifold
$X$, with second order interpolation on any closed complex submanifold
$X_0 \subset X$, then $Y$ admits a dominating spray} [Gro3, FP3].
In [Gro3] the reader can find some further examples of target
manifolds $Y$ for which this holds, but it is not known whether 
it holds for all manifolds.

%
%
%  REMOVING INTERSECTIONS WITH COMPLEX SUBVARIETIES
%
%
\beginsection 4. Removing intersections with complex subvarieties.

Let $X$ and $Y$ be complex manifolds and $A\subset Y$ a closed complex
subvariety of $Y$. Given a \holo\ map $f\colon X \to Y$ we write
$f^{-1}(A)=\{x\in X\colon f(x)\in A\}$ and call it the 
{\it intersection set of $f$ with $A$}. The question is to what
extent is it possible to prescribe the intersection set if $f$ is 
allowed to vary within a homotopy class of maps $X\to Y$. 
More precisely, we consider the following

\smallskip\ni\bf 4.1 Problem. \rm
Suppose that $f^{-1}(A)=X_0\cup X_1$, where $X_0, X_1\subset X$ are 
disjoint complex subvarieties of $X$. When is it possible to remove $X_1$ 
from $f^{-1}(A)$ by homotopy of holomorphic maps $f_t\colon X\to Y$
($t\in [0,1]$) which is fixed on $X_0$ and satisfies $f_0=f$ and 
$f_1^{-1}(A)=X_0$?
\smallskip

In the simplest case when $X=\C$ and $A$ consists of $d$ points in $Y= \CP^1$ 
the answer changes when passing from $d=2$ to $d=3$: 
One can prescribe the pull-back of any two points in $\CP^1$ by a \holo\ map 
$f\colon \C\to\CP^1$ (and there are infinitely many such maps), but when 
$d\ge 3$ the pull-back divisor $f^*A$ completely determines the map $f$. 
Similar situation occurs when $A$ consists of $d$ hyperplanes in general 
position in $Y=\CP^n$: we have flexibility ($=$infinitely many maps) 
for $d \le n+1$ and rigidity (few maps) for $d\ge n+2$.

We say that the Oka principle holds if the existence
of a homotopy of continuous maps $X\to Y$ (which remain holomorphic near 
$X_0$ and remove $X_1$ from the preimage) implies the existence
of a holomorphic homotopy with the required properties.
We break down the problem as follows.

\itemitem{Step 1:} Find a homotopy $f_t\colon X\to Y$ $(0\le t\le 1/2)$ with
$f_0=f$ such that each $f_t$ equals $f$ in an open \nbd\ of $X_0$
in $X$ and $f_{1/2}^{-1}(A)=X_0$. This is a homotopy theoretical problem.

\itemitem{Step 2:} With $f_{1/2}$ as in Step 1,
find a homotopy $f_t\colon X\to Y$  $(1/2\le t \le 1)$ 
such that each $f_t$ is \holo\ near $X_0$ and matches $f$ on $X_0$,
$f_t^{-1}(A)=X_0$ for each $t$, and $f_1$ is \holo\ on $X$.
A solution is given by Theorem 4.2 below when $X$ is Stein and 
$Y\bs A$ is subelliptic.

\itemitem{Step 3:} Deform the combined homotopy $f_t$ $(0\le t\le 1)$ 
from Steps 1 and 2, with fixed $f_0$ and $f_1$, 
to a holomorphic homotopy $\wt f_t\colon X\to Y$ $(t\in [0,1])$ 
such that the resulting two-parameter homotopy is fixed along $X_0$. 
This is possible if $X$ is Stein and $Y$ is subelliptic (Theorem 4.4 below).

\smallskip 
The conclusion is that {\it the Oka principle holds in Problem 4.1 provided that 
$X$ is Stein and the manifolds $Y$ and $Y\bs A$ are subelliptic.} 
Examples of such pairs $A\subset Y$ are algebraic subvarieties 
of codimension at least two (or those with homogeneous complement)
in $\C^n$, $\CP^n$, or in a complex Grassmanian manifold.

We now present the results mentioned above (see [F4, F5] for more).
Thus Theorem 4.2 provides a solution to Step 2 and Theorem 4.4 provides 
a solution to Step 3.

\proclaim 4.2 Theorem: Let $A$ be a closed complex subvariety of a complex
manifold $Y$. Let $X$ be a Stein manifold, $K \subset X$ a compact 
$\cO(X)$-convex subset, and $f\colon X\to Y$ a continuous map which 
is \holo\ in an open set $U_0\subset X$ containing $f^{-1}(A) \cup K$. 
If $Y\bs A$ is subelliptic then for any\ $r\in\N$ there are an open set 
$U\supset f^{-1}(A)\cup K$ and a homotopy $f_t\colon X\to Y$ $(t\in[0,1])$ 
of continuous maps such that $f_0=f$, $f_t$ is \holo\ in $U$ and tangent to 
$f$ to order $r$ along $f_t^{-1}(A)=f^{-1}(A)$ for each $t\in [0,1]$, 
and $f_1$ is \holo\ on $X$.

\ni\ \bf 4.3 Corollary. \sl  
The conclusion of Theorem 4.2 holds in the following cases:
\item{(a)} $Y$ is an affine space $\C^n$, a projective space $\CP^n$
or a complex Grassmanian and $A\subset Y$ is an algebraic subvariety of 
codimension at least two.
\item{(b)} $Y=\CP^n$ and $A$ consists of at most $n+1$ hyperplanes 
in general position.
\item{(c)} A complex Lie group acts transitively on $Y\bs A$. 
\smallskip\rm

In any of these cases $Y\bs A$ is subelliptic by the results stated in 
section three. Note that (b) is a special case of (c).

\proclaim 4.4 Theorem: {\bf (The Oka principle for removing intersections.)} 
Assume that $f\colon X\to Y$ is a \holo\ map of a Stein manifold $X$
to a subelliptic manifold $Y$, $A\subset Y$ is a complex subvariety 
such that $Y\bs A$ is subelliptic, and $f^{-1}(A)=X_0 \cup X_1$ where $X_0$ 
and $X_1$ are unions of connected components of $f^{-1}(A)$ with 
$X_0\cap X_1=\emptyset$. If there exists a homotopy 
$\wt f_t\colon X\to Y$ $(t\in [0,1])$ of continous maps satisfying 
$\wt f_0=f$, $\wt f_1^{-1}(A)=X_0$, and $\wt f_t|_U=f_t|_U$ 
for some open set $U\supset X_0$ and for all $t\in [0,1]$ then for each 
$r\in\N$ there exists a homotopy of \holo\ maps $f_t\colon X\to Y$ such that 
$f=f_0$, $f_1^{-1}(A)=X_0$, and for each $t\in [0,1]$ the map $f_t$ 
agrees to order $r$ with $f$ along $X_0$.

Note that $X_0$ is a union of connected components of $f_t^{-1}(A)$ 
for every $t\in [0,1]$.
In plain language Theorem 4.4 says the following. Suppose that we can remove 
$X_1$ from  $f^{-1}(A)=X_0 \cup X_1$ by a homotopy of continuous maps $X\to Y$ 
which agree with $f$ in a \nbd\ of $X_0$. If $Y$ and $Y\bs A$ are both subelliptic 
then $X_1$ can also be removed from the preimage of $A$ by a homotopy of \holo\ 
maps $X\to Y$ which agree with $f$ to any given order on $X_0$. 

Theorem 4.4 applies if $Y$ is any of the manifolds $\C^n$, $\CP^n$ or a complex 
Grassmanian (these are complex homogeneous and therefore elliptic) and $A\subset Y$ 
is as in Corollary 4.3. When $Y=\C^d$ and $Y\bs A$ is elliptic, Theorem 4.4 coincides 
with Theorem 1.3 in [F5].

%
%  THE OKA PRINCIPLE FOR COMPLETE INTERSECTIONS
%
\smallskip\ni\bf 4.5 Example: The Oka principle for complete intersections. \rm
When $A$ is the origin in $Y=\C^d$ Theorem 4.4 implies the following result of 
Forster and Ramspott [FRa] from 1967. Suppose that a complex subvariety $X_0 \subset X$ 
is a {\it complete intersection} of codimension $d$ in an open set $U\supset X_0$ 
in $X$, given by $d$ functions $f_1,\ldots, f_d \in \cO(U)$ which together
generate the ideal of $X_0$ at each point. If these functions admit continuous 
extensions to $X$ with no additional common zeros then $X_0$ is a (global) 
complete intersection in $X$. The analogous result holds for set-theoretic 
complete intersections.

%
%  SMOOTH VERSUS HOLOMORPHIC COMPLETE INTERSECTIONS
%
\smallskip\ni\bf 4.6 Example: Smooth versus holomorphic complete intersections. \rm
In [F4] it was proved that {\it there exists a three dimensional closed complex 
submanifold $X$ in $\C^5$ which is a smooth (even \ra) complete intersection but which
is not a holomorphic complete intersection.} More precisely, given any compact 
orientable two dimensional surface $M$ of genus $g\ge 2$, there is a three 
dimensional Stein manifold $X$ which is homotopy equivalent to $M$ and whose 
tangent bundle $TX$ is trivial as a real vector bundle but is nontrivial as a 
complex vector bundle over $X$. The image of any \phe\ of $X$ in $\C^5$ (or in $\C^7$) 
is a smooth complete intersection in $\C^5$ (resp.\ $\C^7$) but is not a holomorphic 
complete intersection in any open \nbd\ of $X$ (since its normal bundle is nontrivial 
as a complex vector bundle on $X$). The following problem remains open.

\ni \bf Problem: \sl Let $X\subset \C^n$ be a closed complex submanifold
such that (i) $X$ is a smooth complete intersection in $\C^n$, and 
(ii) its normal bundle $T\C^n|_X /TX$ is trivial as a complex vector bundle
(hence $X$ is a \holo\ complete intersection in an open \nbd\ $U\subset \C^n$). 
Is $X$ a \holo\ complete intersection in $\C^n$~? \rm

%
%  UNAVOIDABLE DISCRETE SETS
%
\smallskip\ni\bf 4.7 Example: Unavoidable discrete sets. \rm 
Theorem 4.4 fails if $Y=\C^n$ and $A$ is any unavoidable discrete subset 
of $\C^n$ (see Example 2.7 above). To see this, write $A=\{p\}\cup A_1$ for some 
$p\in A$. Then $A_1$ is still unavoidable and consequently every entire map 
$F\colon \C^n\to\C^n\bs A_1$ has rank $<n$ at each point. Take $X=\C^n$, 
$f=Id \colon\C^n\to\C^n$, $X_0=\{p\}$ and $X_1=A_1$. The conditions of Theorem 4.4 
are clearly satisfied but its conclusion fails since the rank condition for holomorphic 
maps $F\colon \C^n\to \C^n\bs A_1$ implies that $F^{-1}(p)$ contains no isolated points,
and hence $X_0=\{p\}$ cannot be a connected component of $F^{-1}(p)$.

%
%
% IMMERSIONS AND EMBEDDINGS OF STEIN MANIFOLDS
%
%
\beginsection 5. Embeddings and immersions of Stein manifolds.

In this section we collect results on holomorphic immersions and 
embeddings of Stein manifolds into affine complex spaces. There has been
considerable progress on this subject since the 1990 survey of
Bell an Narasimhan [BN]. 

In 1956 Remmert [Rem] proved that every Stein manifold of dimension 
$n\ge 1$ admits a \phe\ in $\C^{2n+1}$ and a proper \holo\ immersion 
in $\C^{2n}$. Further results were obtained by several authors
(Narasimhan [Na1, Na2], Bishop [Bis], Ramspott [Ram], Forster [Fs1, Fs2], 
Schaft [Sht], and others). The following optimal result, due to Eliashberg 
and Gromov, was announced in 1971 [GE] and proved in 1992 [EG]). An improvement 
of the embedding dimension for odd $n$ is due to Sch\"urmann [Sch] (1997).

\proclaim 5.1 Theorem: Every Stein manifold $X$ of dimension $n>1$ admits 
a proper holomorphic embedding in $\C^{[3n/2]+1}$ and a proper \holo\ 
immersion in $\C^{[3n+1/2]}$. Also there exists a (not necessarily proper)
holomorphic immersion $X\to \C^{[3n/2]}$.

Sch\"urmann [Sch] also proved an optimal embedding theorem 
for Stein spaces with singularities which have uniformly bounded local embedding 
dimension. Recently J.\ Prezelj [Pr2] constructed {\it proper weakly holomorphic 
embeddings of Stein spaces with isolated singularities in Euclidean spaces of 
minimal dimension.} 
These results strongly depend on Lefschetz's theorem to the effect
that a Stein manifold is homotopically equivalent to a CW-complex  
of real dimension at most $\dim_{\C} X$ [AF]. The following example of 
Forster [Fs1] shows that the dimensions in Theorem 5.1 cannot be lowered. 

%
% Forster's example
%
\smallskip\ni \bf 5.2 Example. \rm Let
$$ \eqalign{ Y &= \{[x\colon y\colon z]\in \CP^2\colon X^2+y^2+z^2\ne 0\} \cr
             X&= \cases { Y^m, & if $n=2m$;\cr
                          Y^m\times\C, & if $n=2m+1$. \cr}
           }               
$$
Clearly $X$ is a Stein manifold of dimension $n$. Forster proved 
[Fo1, Proposition 3] that the Stiefel-Whitney class $w_{2m}(TX)$ 
is the nonzero element of the group $H^{2m}(X;Z_2) = H^{2m}((\RP^2)^m; Z_2)=Z_2$
and consequently the Chern class $c_m(TX)$ is the nonzero element 
of $H^{2m}(X;Z)=Z_2$. It follows that $X$ does not embed in $\C^{[3n/2]}$ 
and does not immerse in $\C^{[3n/2]-1}$ (see [Hu, p.\ 263]). 
\endpr

The proof of Theorem 5.1 in [EG] and [Sch] relies on the elimination 
of singularities method and on the Oka principle for sections of certain 
submersions onto Stein manifolds (Example 3.9 above). 
We describe the main idea. One begins by choosing a generic almost proper holomorphic map 
$b\colon X^n\to\C^n$ constructed by Bishop [Bis]. This means that the $b$-preimage of 
any compact set in $\C^n$ has (at most countably many) compact connected components.
We then try to find a map $g\colon X\to \C^q$ which `desingularizes $b$' in the 
sense that $f=(b,g)\colon X\to\C^{n+q}$ is a proper holomorphic embedding
(resp.\ immersion). Properness is easily achieved by choosing $g$ sufficiently 
large on a certain sequence of compact sets in $X$ (here we need that $b$ is almost proper). 

To insure that $f=(b,g)$ is an immersion we must choose $g$ such that its differential 
$dg_x$ is nondegenerate on the kernel of $db_x$ at each $x\in X$. To obtain
injectivity we must choose $g$ to separate points on the fibers of $b$. 
Both requirements can be satisfied if $q\ge [n/2]+1$
and this number is determined by topological restrictions. (The immersion
condition  requires $q\ge [n/2]$.) One proves this by a finite induction. 
We stratify $X$ by a descending finite chain of closed complex subvarieties 
$X=X_0\supset X_1\supset X_2\ldots\supset X_m=\emptyset$ such that the kernel of 
$db_x$ has constant dimension on each stratum $S_k=X_k\bs X_{k+1}$ (which is chosen
to be nonsingular), and the number of distinct points in $b^{-1}(x)$ is 
constant for $x\in S_k$. (In the actual proof we must replace $X$ 
by a suitable subset $B \subset X$ which is mapped by $b$ properly onto a 
bounded domain in $\C^n$; in the end we perform an induction by increasing $B$ 
to $X$.) Furthermore, once we have a map $g_k\colon X\to\C^q$ satisfying these 
conditions along $X_k$, we choose $g_{k-1}\colon X\to\C^q$ such that it satsfies
both condition on the next stratum $S_{k-1}$ and agrees with $g_k$ to second 
order along $X_k$ (so that $g_{k-1}$ does not destroy what $g_k$ has achieved). 
A suitable $g_{k-1}$ is obtained by the Oka principle (Section 3) provided 
there are no topological obstructions, and this is so when $q\ge [n/2]+1$.

Although Theorem 5.1 gives the optimal result for the entire collection 
of $n$-dimensional Stein manifolds, the method does not give a better result
for `simple' Stein manifolds which are expected to embed in lower dimensional
space. For instance, it is not known what is the minimal proper embedding 
dimension of the polydisc or the ball in $\C^n$. Globevnik proved 
by a different method (using shear automorphisms of $\C^n$) that there are 
arbitrarily small perturbations of the polydisc in $\C^n$ which embed 
in $\C^{n+1}$ [Gl2].

\proclaim 5.3 Question: 
What is the proper holomorhic embedding dimension of the ball?
The polydisc? A general convex domains in $\C^n$? 
How does it depend on the topology and geometry of the domain?

We now consider the existence of {\it relative embeddings}.
The following result was proved in [ABT] following the method
of Narasimhan [Na2]. 

\proclaim 5.4 Theorem: Suppose that $X$ is a Stein manifold of dimension $n$, 
$Y\subset X$ is a closed complex submanifold in $X$ and $f\colon Y\to\C^N$ 
is a \phe. If $N\ge 2n+1$ there exists a \phe\ $\wt f\colon X\to\C^N$ 
extending $f$.

It is not known whether Theorem 5.4 is valid for $N=2n$, but it is  false 
for $N \le 2n-1$ by Corollary 5.7 below which follows from the following
interpolation results for holomorphic embeddings from [BFn] and [F3].

%
%
%  EMBEDDINGS WITH INTERPOLATION
%
%

\proclaim 5.5 Theorem:
Let $\Sigma$ be a discrete subset of $\C^N$ for some $N>1$.
If a Stein manifold $X$ admits a proper holomorphic embedding
$f_0\colon X\to\C^N$ then $X$ also admits an embedding $f\colon X\to \C^N$ 
whose image $f(X)$ contains $\Sigma$. In addition we may choose $f$ such 
that for every entire map $\psi\colon \C^d\to\C^N$ whose rank
equals $d=N-\dim X$ at most points of $\C^d$ the set $\psi(\C^d)\cap f(X)$ 
is infinite. If $d=1$, we may insure that $\C^N\bs f(X)$ is Kobayashi hyperbolic.

\proclaim 5.6 Corollary:
Let $n,d\ge 1,\ N=n+d$. There exists a \phe\ $f\colon \C^n\to\C^N$
such that every entire map $\psi\colon \C^d\to\C^N$ of rank $d$
intersects $f(\C^n)$ at infinitely many points. For $d=1$ 
we may choose $f$ such that $\C^{n+1}\bs f(\C^n)$ 
is Kobayashi hyperbolic.

The proofs in [BFn] and [F3] use results on holomorphic automorphisms of 
$\C^N$ obtained in [And, AL, FRo]. The first result in this direction 
[FGR] was that {\it there exist holomorphically embedded complex lines in $\C^2$ 
which are not equivalent to the standard embedding $\C\to \C\times\{0\}\subset\C^2$ by 
automorphisms of $\C^2$}. This is in strong contrast to the situation 
for algebraic (polynomial) embeddings $\C\to\C^2$ which are all equivalent to
the standard embedding by polynomial automorphisms of $\C^2$ according to
Abhyankar and Moh [AM].  Using such `twisted' holomorphic embeddings
of $\C$ in $\C^2$ Derksen and Kutzschebauch constructed 
nonlinearizable periodic holomorphic automorphisms of $\C^4$ [DK].

\proclaim 5.7 Corollary: For every $n\ge 2$ there exists a \phe\
$f\colon \C^{n-1}\to \C^{2n-1}$ which does not admit an injective 
holomorphic extension $\wt f\colon \C^n\to \C^{2n-1}$.

\demo Proof: Choose $f\colon \C^{n-1}\to \C^{2n-1}$ as in Corollary 5.6 such
that the range of any entire map $\C^{n}\to \C^{2n-1}$ of generic rank $n$ 
intersects $f(\C^{n-1})$. If $\wt f\colon \C^n\to\C^{2n-1}$ is an 
injective \holo\ extension of $f\colon \C^{n-1}\times\{0\}\to\C^{2n-1}$ then 
$$ 
	\psi(z)=\psi(z_1,\ldots,z_n)= \wt f\bigl(z_1,\ldots,z_{n-1},e^{z_n}\bigr) 
	\qquad (z\in\C^n)
$$
is an entire map which has rank $n$ at a generic point of $\C^n$ and whose image 
misses $f(\C^{n-1}) \subset \C^{2n-1}$ in contradiction to the assumption on $f$.
\endpr

In Theorem 5.5 the image $f(X)\subset \C^N$ contains a given discrete 
set $\{p_j\} \subset \C^N$, but we don't specify the points in $X$ which 
correspond to the points $p_j$ under the embedding $f$. The following
more precise interpolation theorem was proved recently by J.\ Prezelj 
[Pr1]. 

%
% Jasna's theorem
%
\proclaim 5.8 Theorem: Let $X$ be a Stein manifold of dimension $n\ge 1$.
Define $q(n)= \min\{ [{n+1\over 2}]+1, 3\}$. Then  for any
$N\ge n+q(n)$ and any pair of discrete sets $\{a_k\}\subset X$,
$\{b_k\}\subset \C^N$ there exists a \phe\ $f\colon X\to\C^N$
satisfying $f(a_k)=b_k$ for every $k=1,2,3,\ldots$. The analogous 
conclusion holds for proper \holo\ immersions $X\to\C^N$ when 
$N\ge [{3n+1 \over 2}]$.

Comparing with Theorem 5.1 we see that the embedding dimension is minimal 
for even $n$ and is off by at most one for odd $n$.  Prezelj's proof in 
[Pr1] uses an improved version of the scheme from [EG] and [Sch]. 
By entirely different methods (using holomorphic automorphisms)
J.\ Globevnik proved that the conclusion of Theorem 5.8 also holds for 
proper holomorphic embeddings of the unit disc in $\C^2$ [Gl3]. A Carleman type 
embedding theorem (approximating a given smooth proper embedding $\R\to\C^n$ 
in the fine $\cC^\infty$ topology on $\R$ by proper \holo\ embeddings $\C\to\C^n$) 
was proved in [BF].

%
%
%  THE HOMOTOPY PRINCIPLE FOR IMMERSIONS
%
%

It is not known whether {\it proper} holomorphic immersions or embeddings
of Stein manifolds satisfy the Oka principle. However, non-proper holomorphic 
immersions of Stein manifolds do satisfy the following Oka principle 
(Gromov and Eliashberg [GE]; see also section 2.1.5.\ in [Gro2]).

%
%  Gromov and Eliashberg
%
\proclaim 5.9 Theorem:  If the cotangent bundle $T^*X$ of a Stein manifold
is generated by $q$ differential $(1,0)$-forms $\theta_1,\ldots,\theta_q$
for some $q>\dim X$ then there exists a holomorphic immersion $X\to\C^q$. 
More precisely, every such $q$-tuple $(\theta_1,\ldots,\theta_q)$ 
can be changed by a homotopy (through $q$-tuples generating $T^*X$) 
to a $q$-tuple $(df_1,\ldots,df_q)$ where $f=(f_1,\ldots,f_q)\colon X\to\C^q$ 
is a \holo\ immersion.

The idea of the proof is the following. By the Oka-Grauert principle we may 
assume that $\theta_j$ are holomorphic 1-forms.
In the first step one of the forms, say $\theta_q$, is replaced by the 
differential $df_q$ of a  holomorphic function on $X$ such that 
$\theta_1,\ldots,\theta_{q-1},df_q$ still generate $T^*X$. Since 
$q>\dim X$, we may assume that the forms $\theta_1,\ldots,\theta_{q-1}$
already generate $T^*X$ outside a proper complex subvariety
$\Sigma \subset X$, and $f_q$ must satisfy an essentially algebraic 
condition on its jet along $\Sigma$. Once $f_q$ has been chosen one proceeds
in the same way and replaces $\theta_{q-1}$ with an exact differential.
In finitely many steps all forms are replaced with differentials.
The technical details of the proof are considerable.

%
%
%  EMBEDDINGS OF OPEN RIEMANN SURFACES IN THE AFFINE PLANE
%
%
\beginsection 6. Embeddings of open Riemann surfaces in the affine plane.

In this section we describe the state of knowledge on the following problem.

\proclaim 6.1 Problem: Does every open Riemann surface admit a proper
holomorphic embedding in $\C^2$? Is the algebra of global holomorphic functions
on such a surface always doubly generated?

Open Riemann surfaces are precisely Stein manifolds of dimension one,
and in view of Theorem 5.1 (on embedding $n$-dimensional Stein manifolds in 
$\C^{[3n/2]+1}$ for $n>1$) one might expect that they embed in $\C^2$. 
(For comparison we recall that every compact Riemann surface embeds in 
$\CP^3$ but most of them don't embed in $\CP^2$ [FK].) The proof 
of Theorem 5.1 only gives embeddings into $\C^3$, the reason being that 
for embeddings $X\hra\C^2$ it runs into a hyperbolicity obstruction 
(Example 2.6 in Section 2).  The main difficulty is 
to find injective holomorphic maps to $\C^2$;
this is essentially equivalent to the algebra of holomorphic functions 
being doubly generated. Here are some Riemann surfaces which are known 
to embed in $\C^2$:
\item{--} the disc $U=\{z\in\C\colon |z|<1\}$ (Kasahara and Nishino [Ste]);
\item{--} annuli $\{1<|z|<r\}$ (Laufer [Lau]),
\item{--} punctured disc $U\bs \{0\}$ (Alexander [Ale]),
\item{--} all finitely connected planar domains $\Omega\subset \C$ different from
$\C$ whose boundary contains no isolated points (Globevnik and Stens\o nes [GS]).

\smallskip
We now consider the embedding problem for {\it bordered Riemann surfaces}. 
Let $\cR$ be a compact, orientable, smooth real surfaces whose boundary 
$b\cR=\cup_{j=1}^m C_j$ consists of finitely many curves and no isolated points. 
Such $\cR$ is a sphere with $g$ handles ($g$ is the {\it geometric genus} 
of $\cR$) and $m\ge 1$ holes (removed discs). A {\it complex structure} on $\cR$ is 
determined by a real endomorphism $J$ of the 
tangent bundle $T\cR$ satisying $J^2=-Id$ (Gauss-Ahlfors-Bers). 
Without loss of generality we may assume that $J$ is H\"older continuous of 
class $\cC^\a(\cR)$ for a fixed $\a\in(0,1)$. A differentiable function 
$f\colon \cR\to\C$ is $J$-holomorphic if $df\circ J=i df$ where $i=\sqrt{-1}$. 
Two complex structures $J_0$ and $J_1$ are {\it equivalent} if there exists 
a diffeomorphism $\phi\colon \cR\to\cR$ of class $\cC^{1,\a}(\cR)$ satisfying 
$d\phi\circ J_0=J_1\circ d\phi$. The set of equivalence classes of complex 
structures on $\cR$ is the {\it moduli space} $\cM(\cR)$ of Riemann surface 
structures on $\cR$. The following result from [\v CF] shows that 
{\it there are no topological obstructions for embedding finite bordered
Riemann surfaces in $\C^2$}.

\proclaim 6.2 Theorem: For every smooth bordered surface $\cR$ there exists a 
nonempty open set  $\Omega\subset \cM(\cR)$ such that for every 
complex structure $J$ on $\cR$ with $[J]\in \Omega$ the open Riemann surface 
$\IR=\cR\bs b\cR$ admits a proper $J$-holomorphic embedding in $\C^2$.

Theorem 6.2 follows from the following result in [\v CF] which sems to contain
all known results on embeddings in $\C^2$ except for the punctured disc.
(For planar domains see [GS] and [\v CG]). Let $U=\{z\in\C\colon |z|<1\}$.

\proclaim 6.3 Theorem:  Let $(\cR,J)$ be a finite bordered Riemann surface 
of genus $g$ with $m$ boundary components, where $J$ is of class $\cC^\a(\cR)$
for some $\a\in (0,1)$.  Assume that there exists an 
injective immersion $f=(f_1,f_2) \colon \cR \to \bar U\times \C$ of class $\cC^2$ 
which is $J$-\holo\ in $\IR$, $|f_1|=1$ on $b\cR$, 
and the generic fiber of $f_1$ contains at least $2g+m-1$ points. 
Then $\IR$ admits a proper $J$-holomorphic embedding in $\C^2$. 
Furthermore, for every complex structure $\wt J$ sufficiently $\cC^\a$ close 
to $J$ the surface $\IR$ also admits a proper $\wt J$-holomorphic embedding in $\C^2$.

\smallskip\ni\bf 6.4 Corollary. \sl
The following Riemann surfaces admit a proper \holo\ embedding in $\C^2$:
\item{(i)} finitely connected domains in $\C$ without isolated 
boundary points, 
\item{(ii)} every complex torus with one hole, 
\item{(iii)} every bordered Riemann surface whose double is hyperelliptic.
\smallskip\rm

The proof of Theorems 6.2 and 6.3 in [\v CF] is based partly on the method 
developed by Globevnik and Sten\o nes [GS] who proved the result for planar domains.
In this case the conditions in Theorem 6.3 are satisfied if we take $g(z)=z$ 
and $f$ an inner function of degree $\ge m-1$, where $m$ is the number 
of boundary components of the domain.

A hyperelliptic (compact) Riemann surface $X$ is the normalization of a curve 
in $\CP^2$ given by 
$$ 
	y^2=\prod_{j=0}^{g} (x-\a_j)(1-\bar \a_j x) \eqno(6.1)
$$
for distinct points $\a_j\in U$ ($0\le j\le g$), where $g=g(X)$ is the genus 
of $X$. If $X$ is the double of a bordered Riemann surface $\cR$ then 
$g(X)=2g(\cR)+m-1$ where $m$ is the number of boundary curves of $\cR$ 
(which equals either $1$ or $2$ in this case), and the representation (6.1) 
can be chosen such that $\cR=\{(x,y)\in X \colon |x|\le 1\}$. 
The pair of functions 
$$
	f_1 = y/\prod_{j=0}^{g} (1-\bar\a_j x), \quad f_2= x
$$ 
provides an embedding $f=(f_1,f_2)\colon \cR \to \bar U^2$ which maps $b\cR$ 
into the torus $(bU)^2$. Clearly $g$ has multiplicity two. From 
$f_1^2=\prod_{j=0}^{g} (f_2 - \a_j)/(1-\bar\a_j f_2)$ which follows 
from (6.1) we see that $f_1$ has multiplicity $g+1=2g(\cR)+m$ 
and hence Theorem 6.3 applies. Sikorav gave a slightly different proof 
for embedding of complex tori with one hole (unpublished); 
these are all hyperelliptic.

The proof of Theorem 6.3 goes as follows.
Let $P=(2U)\times (RU) \subset \C^2$ for some $R>\sup |f_2|$.  
Globevnik and Stens\o nes [GS, Gl1] found arbitrarily small smooth perturbations 
$S\subset \C^2$ of the cylinder $S_0=bU\times \C$ such that the connected component 
$\Omega_S$ of $P\bs S$ containing the origin is of the form 
$\Omega_S=\wt \Omega_S \cap P$ for some Fatou-Bieberbach domain 
$\wt \Omega_S \subset\C^2$. Let $\phi_S\colon \wt\Omega_S\to\C^2$ be a biholomorphic 
(Fatou-Bieberbach) map. If $f^S=(f^S_1,f^S_2)\colon \cR\to \bar \Omega_S$ is a 
continuous map which maps $\IR$ holomorphically into $\Omega_S$ and maps $b\cR$ 
into $S\cap P$  then clearly $\phi_S \circ f^S\colon \IR\to\C^2$ is a \phe. 
A map $f^S$ with these properties is obtained from the initial map $f$ 
satisfying the hypothesis of Theorem 6.3  by solving a Riemann-Hilbert boundary 
value problem on $\cR$ (in fact we only need to perturb the first component of $f$).

\smallskip\ni\bf 6.5 Open problems. \sl
\item{(a)} Does every finite Riemann surface with boundary consisting of 
finitely many closed curves and isolated points embed in $\C^2$?
\item{(b)} Does every planar domain with finitely many punctures 
embed in $\C^2$?
\smallskip\rm

\beginsection 7. Noncritical holomorphic functions and submersions.

In 1967 Gunning and Narasimhan [GN] proved that every open Riemann surface admits 
a \holo\ function without critical points, thus giving an affirmative answer 
to a long standing question. Open Riemann surfaces are precisely Stein 
manifolds of complex dimension one. In the recent paper [F8] the result 
of [GN] was extended to Stein manifolds of any dimension.

\proclaim 7.1 Theorem: {\rm (Theorem I in [F8].)}
Every Stein manifold admits a \holo\ function without critical points.
More precisely, an $n$-dimensional Stein manifold admits 
$[{n+1\over 2}]$ holomorphic functions with pointwise independent 
differentials and this number is maximal for every $n$.

In fact the manifold $X$ in Example 5.2 above does not admit 
$[{n+1\over 2}]+1$ functions with independent differentials,
the reason being that $c_m(X)\ne 0$ for  $m=[n/2]$
(Proposition 2.11 in [F8]).

Furthermore it is proved in [F8] (Theorem 2.1 and  Corollary 2.2)
that {\it for any discrete subset $P$ in a Stein manifold $X$ there exists 
a \holo\ function $f\in\cO(X)$ whose critical set equals $P$, and one can prescribe 
the finite order jet of $f$ at each point of $P$}. Any noncritical holomorphic function 
on a closed complex submanifold $X_0$ in a Stein manifold $X$ extends to a 
noncritical holomorphic function on $X$, and hence  $X$ admits a 
nonsingular hypersurface foliation transverse to the submanifold $X_0$.

The main result of [F8] is a {\it homotopy principle for holomorphic 
submersions of Stein manifolds to Euclidean spaces of lower dimension.}
A \holo\ map $f\colon X\to\C^q$ is a {\it submersion} if 
$df_x\colon T_x X\to T_{f(x)}\C^q \simeq \C^q$ is surjective for every $x\in X$,
and hence its differential induces a surjective complex vector bundle map 
of the tangent bundle $TX$  onto the trivial bundle $X\times \C^q$.  
For $q<\dim X$ this necessary condition for the existence of a submersion 
$X\to\C^q$ is also sufficient. The following is Theorem II in [F8].

\proclaim 7.2 Theorem:
If $X$ is a Stein manifold and $1\le q < \dim X$ then every surjective complex 
vector bundle map $TX\to X\times\C^q$ is homotopic (in the space of surjective
vector bundle maps) to the differential of a \holo\ submersion $X\to\C^q$.
If $f_0,f_1\colon X\to\C^q$ are \holo\ submersions whose 
differentials are homotopic through a family of surjective complex vector bundle 
maps $\theta_t\colon TX\to X\times\C^q$ $(t\in [0,1])$ then there exists a
regular homotopy of \holo\ submersions $f_t\colon X\to\C^q$ $(t\in [0,1])$ 
connecting $f_0$ to $f_1$.

For a more precise result (with approcimation and interpolation) see Theorem 2.5 in [F8]. 
It is not known whether the same conclusion holds for $q=\dim X$, 
except on open Riemann surfaces where it was proved in [GN]. 
The corresponding homotopy principle for submersions of smooth open manifolds 
to $\R^q$ is due to Phillips [Ph1, Ph3] (see also [Gr2]). For an extension to 
smooth (or real-analytic) foliations see [Ph2] and [Th1, Th2].

In [F8] the reader can find numerous applications to the existence 
of nonsingular \holo\ foliations on Stein manifolds. For instance
we have the following (Corollary 2.9 and Theorem 7.1 in [F8]).

\proclaim 7.3 Corollary: 
Let $X$ be a Stein manifold of dimension $n$ and $E\subset TX$ a complex 
subbundle of rank $k\ge 1$. If the quotient bundle $TX/E$ is trivial or
admits locally constant transition functions then $E$ is homotopic 
(through complex rank $k$ subbundles of $TX$) to the tangent bundle of a nonsingular
holomorphic foliation of $X$. If $E$ is holomorphic then the homotopy 
may be chosen through \holo\ subbundles.

When $TX/E$ is trivial (of rank $q=n-k$) the foliation 
in Corollary 7.3 is given by a \holo\ submersion $X\to\C^q$. 

To any complex subbundle $E \subset TX$ we associate its 
polar $\Theta=E^\perp \subset T^*X$ with fibers
$\Theta_x=\{\omega\in T_x^*X\colon \omega(v)=0\ {\rm for\ all\ } v\in E_x\}$.
Then $\Theta =(TX/E)^*$ and the first part of Corollary 7.3 can 
be equivalently expressed as follows:
{\it Any trivial complex subbundle $\Theta\subset T^* X$ of rank $q<\dim X$
is homotopic to a subbundle generated by independent holomorphic differentials
$df_1,\ldots, df_q$.} 

By the Lefschetz theorem [AF] a Stein manifold $X$ of complex dimension $n$ 
is homotopic to a CW-complex of real dimension at most $n$. Elementary 
homotopy theory implies that every complex vector bundle of rank $m\ge [n/2]+1$ 
on $X$ admits a nonvanishing section. 
Applying this inductively we see that  $T^*X$ admits a trivial 
complex subbundle of rank $[(n+1)/2]$. Hence Corollary 7.3 implies

\proclaim 7.4 Corollary: Every $n$-dimensional Stein manifold $X$ admits
nonsingular \holo\ foliations of any dimension $k \ge [n/2]$. 
If $X$ is holomorphically parallelizable, it admits a \holo\ submersion 
$X\to \C^{n-1}$ and holomorphic foliations of any dimension 
$1,2,\ldots,n-1$.

It is conjectured that every parallelizable Stein manifold of dimension 
$n$ admits a \holo\ submersion ($=$immersion) in $\C^n$ (see 7.5 below). 
The construction of submersions $X\to\C^q$ in [F8] breaks down for $q=\dim X$
due to a {\it Picard type obstruction} in the approximation 
problem for submersions on Euclidean spaces. 
For example, the complex $n$-sphere 
$\Sigma^n= \{z\in \C^{n+1}\colon \sum_{j=0}^n z_j^2 = 1  \}$
is (algebraically) parallelizable for any $n\in \N$,
and an explicit holomorphic immersion $\S^n\to\C^n$ was found 
by J.\ J.\ Loeb (see [BN, p.18]). The complement $\CP^2\bs C$ of any 
smooth cubic curve $C\subset \CP^2$ is also algebraically parallelizable, 
but is is not known whether they all immerse in $\C^2$.

The construction of noncritical \holo\ functions and submersions in [F8] 
depends on three main ingredients. The first is a method for approximating noncritical 
\holo\ functions on polynomially convex subsets of $\C^n$ by entire noncritical 
functions. This uses the theory of \holo\ automorphisms of affine complex spaces
from [And, AL, FRo]. A similar method is used for submersions $\C^n\to \C^q$ with 
$1<q<n$, but with a weaker conclusion.  

The second crucial new ingredient is a tool for 
patching holomorphic maps preserving the maximal rank condition. 
If $A,B\subset X$ is a  {\it Cartan pair} in a complex manifold $X$ 
then every biholomorphic map $\gamma$ sufficiently uniformly close to the 
identity in a \nbd\ of $A\cap B$ admits a decomposition $\g=\b\circ\a^{-1}$, 
where $\a$ (resp.\ $\b$) is a biholomorphic map close to 
the identity in a \nbd\ of $A$ (resp.\ of $B$). If $\g$ preserves 
a holomorphic foliation of $X$ then $\a$ and $\b$ can be chosen to preserve 
it as well (Theorem 4.1 in [F8]. 
The map $\g$ arises as a `transition map' between a pair of
\holo\ submersions $f$ resp.\ $g$ defined in a \nbd\ of $A$ (resp.\ of $B$),
satisfying $f=g\circ \g$ near $A\cap B$. 
From $\g=\b\circ\a^{-1}$ we obtain $f\circ \a=g\circ \b$ which gives 
a submersion in a \nbd\ of $A \cup B$.

The construction of \holo\ submersions $X\to\C^q$ for $q>1$ in [F8] 
also uses a new device for crossing the critical points of a \spsh\ exhaustion 
function $\rho\colon X\to\R$. (This is not needed for $q=1$.) 
Let $p\in X$ be a Morse critical point of $\rho$. Using a special case 
of Gromov's {\it convex integration lemma} one obtains a smooth maximal 
rank extension of the map $f$ from a sublevel set $\{\rho\le c\} \subset X$
(where $c<\rho(p)$ is close to $\rho(p)$) to a \tr\ handle $E\subset X$ attached 
to $\{\rho \le c\}$ which describes the change of topology at $p$. 
This forms an intrinsic link between the smooth and holomorphic submersions,
the main point being that the complex submersion differential relation is
{\it ample in coordinate directions} on any totally real submanifold.
One then constructs a smooth family of increasing \spsc\ neigborhoods of 
$\{\rho\le c\} \cup E$ such that the largest one contains the sublevel set 
$\{\rho<c'\}$ for some $c'>\rho(p)$. This reduces the extension problem 
to the noncritical case for a different \spsh\ function and lets us pass the critical 
level of $\rho$ at $p$.

\medskip\ni \bf 7.5 Open problems. 	\rm
\item{1.} Let $X$ be a parallelizable Stein manifold of dimension $n>1$.
Does there exist a \holo\ immersion ($=$submersion) $X^n\to\C^n$?
For $n=1$ the affirmative answer was given by Gunning and Narasimhan [GN] in 1967. 
Not much is known for $n>1$ (see e.g.\ [Na3] and [BN]). A positive answer would 
follow from the positive answer to any of the following problems.

\item{2.} 
Let $B$ be an open convex set in $\C^n$. Is every \holo\ immersion 
($=$submersion) $B\to \C^n$ a limit of entire immersions $\C^n\to \C^n$,
uniformly on compacts in $B$~?  The same question may be asked for maps 
with {\it constant Jacobian}; compare with the {\it Jacobian problem} for 
polynomial maps [BN, p.\ 21]. 

\item{3.} Let $f\colon X^n\to\C^{n-1}$ be a holomorphic submersion
and let $\cF$ denote the foliation $\{f=c\}$, $c\in \C^{n-1}$. Assuming that the 
tangent bundle $T\cF$ is trivial, find a $g\in \cO(X)$ which is noncritical 
on every leaf of $\cF$. (The map $(f,g)\colon X\to\C^n$ is then locally biholomorphic.)

\medskip
The methods used in [F8] are highly transcendental and it would be interesting 
to see how much of this theory survives in the algebraic category.

\medskip\ni\bf  References. \rm
\medskip

\ii{[AM]} S.\ S.\ Abhyankar, T.\ Moh: Embeddings of the line in the plane.
J.\ Reine Angew.\ Math.\ {\bf 276}, 148--166 (1975).

\ii{[ABT]}
F.\ Acquistapace, F.\ Broglia, A.\ Tognoli:
A relative embedding theorem for Stein spaces. 
Ann.\ Scuola Norm.\ Sup.\ Pisa Cl.\ Sci.\ (4) {\bf 2}, 507--522 (1975). 

\ii{[Ale]}  H.\ Alexander: 
Explicit imbedding of the (punctured) disc into $\C^2$. 
Comment.\ Math.\ Helv.\ {\bf 52}, 539--544 (1977).

\ii{[And]} E.\ Anders\'en: 
Volume-preserving automorphisms of $\C^n$.
Complex Variables {\bf 14}, 223-235 (1990).

\ii{[AL]} E.\ Anders\'en, L.\ Lempert:
On the group of holomorphic automorphisms of $\C^n$.
Inventiones Math.\ {\bf 110}, 371--388 (1992).

\ii{[AF]} A.\ Andreotti, T.\ Frankel:
The Lefschetz theorem on hyperplane sections.
Ann.\ Math.\ {\bf 69}, 713--717 (1959).

\ii{[Au]} M.\ Audin: Fibr\'es normaux d'immersions en dimension double, 
points doubles d'immersions lagragiennes et plongements totalement r\'eels. 
Comment.\ Math.\ Helv.\ {\bf 63}, 593--632 (1988). 

\ii{[BN]} S.\ Bell and R.\ Narasimhan:
Proper holomorphic mappings of complex spaces.
Several complex variables, VI, 1--38, 
Encyclopaedia Math.\ Sci.\ {\bf 69}, Springer, Berlin, 1990.

\ii{[Bis]} E.\ Bishop: Mappings of partially analytic spaces. 
Amer.\ J.\ Math.\ {\bf 83}, 209--242 (1961).

\ii{[BFn]} G.\ Buzzard, J.-E.\ Forn\ae ss:
An embedding of $\C$ into $\C^2$ with hyperbolic complement.
Math.\ Ann.\ {\bf 306}, 539--546 (1996).

\ii{[BF]} G.\ T.\ Buzzard, F.\ Forstneri\v c:
A Carleman type theorem for proper holomorphic embeddings. 
Ark.\ Mat.\ {\bf  35}, 157--169 (1997).

\ii{[Ca]} H.\ Cartan: Espaces fibr\'es analytiques.
Symposium Internat.\ de topo\-lo\-gia algebraica, Mexico, 97--121 (1958).
(Also in Oeuvres {\bf 2}, Springer, New York, 1979.)

\ii{[\v CF]} M.\ \v Cerne, F.\ Forstneri\v c: 
Embedding some bordered Riemann surfaces in the affine plane.
Math.\ Res.\ Lett.\ {\bf 9} (2002). 

\ii{[\v CG]} M.\ \v Cerne, J.\ Globevnik:
On holomorphic embedding of planar domains into $\C^2$.
J.\ d'Analyse Math.\ {\bf 8}, 269--282 (2000).

\ii{[DK]} H.\ Derksen, F.\ Kutzschebauch:
Nonlinearizable holomorphic group actions.
Math.\ Ann.\ {\bf 311}, 41--53 (1998).

\ii{[El]} Y.\ Eliashberg: 
Topological characterization of Stein manifolds of dimension $>2$. 
Internat.\ J.\ Math. {\bf 1}, 29--46  (1990). 

\ii{[EG]} Y.\ Eliashberg, M.\ Gromov: Embeddings of Stein manifolds.
Ann.\ Math.\ {\bf 136}, 123--135 (1992).

\ii{[EM]} Y.\ Eliashberg, N.\ Mishachev: 
Introduction to the $h$-principle. 
Graduate Studies in Mathematics {\bf 48},  
Amer.\ Math.\ Soc., Providence, 2002. 

\ii{[FK]} H.\ M.\ Farkas, I.\ Kra: Riemann surfaces. 
Second ed. Graduate Texts in Mathematics {\bf 71}, 
Springer, New York, 1992.

\ii{[Fe]} S.\ Feit: $k$-mersions of manifolds.
Acta Math.\ {\bf 122} 173--195 (1969).

\ii{[Fs1]} O.\ Forster: Plongements des vari\'et\'es de Stein.
Comm.\ Math.\ Helv.\ {\bf 45}, 170--184 (1970).

\ii{[Fs2]} O.\ Forster: 
Topologische Methoden in der Theorie Steinscher R\"aume.
(Proc.\ ICM--Nice 1970, 613--618), Nice, 1971.

\ii{[FRa]}  O.\ Forster and K.\ J.\ Ramspott:
Analytische Modulgarben und Endromisb\"undel.
Invent.\ Math.\ {\bf 2}, 145--170 (1966).

\ii{[F1]} F.\ Forstneri\v c:
On totally real embeddings into $\C^n$.
Expo.\ Math.\ {\bf 4}, 243--255 (1986). 

\ii{[F2]} F.\ Forstneri\v c:
Complex tangents of real surfaces in complex surfaces.
Duke Math.\ J.\ {\bf 67}, 353--376  (1992).

\ii{[F3]} F.\ Forstneri\v c:
Interpolation by holomorphic automorphisms and embeddings in $\C^n$.
J.\ Geom.\ Anal.\ {\bf 9}, 93-118 (1999).

\ii{[F4]} F.\ Forstneri\v c: 
On complete intersections. 
Ann.\ Inst.\ Fourier {\bf 51}, 497--512  (2001).

\ii{[F5]} F.\ Forstneri\v c:
The Oka principle for sections of subelliptic submersions.
Math.\ Z.\ {\bf 241}, 527--551 (2002).

\ii{[F6]} F.\ Forstneri\v c:
The Oka principle for sections of ramified mappings.
Forum Math., to appear.

\ii{[F7]} F.\ Forstneri\v c: 
Stein domains in complex surfaces. 
J.\ Geom.\ Anal., to appear.

\ii{[F8]} F.\ Forstneri\v c: 
Noncritical holomorphic functions on Stein manifolds.
Acta Math., to appear.

\ii{[FP1]} F.\ Forstneri\v c and J.\ Prezelj:
Oka's principle for holomorphic fiber bundles with sprays
(with J.\ Prezelj). Math.\ Ann.\ {\bf 317}, 117-154  (2000).

\ii{[FP2]} F.\ Forstneri\v c and J.\ Prezelj:
Oka's principle for holomorphic submersions with sprays.
Math.\ Ann. {\bf 322}, 633-666  (2002).

\ii{[FP3]} F.\ Forstneri\v c and J.\ Prezelj:
Extending holomorphic sections from complex subvarieties.
Math.\ Z.\ {\bf 236}, 43--68  (2001).

\ii{[FGR]} F.\ Forstneri\v c, J.\ Globevnik and  J.-P.\ Rosay:
Non straightenable complex lines in $\C^2$.
Ark.\ Mat.\ {\bf 34}, 97--101  (1996).

\ii{[FRo]} F.\ Forstneric, J.-P.\ Rosay:
Approximation of biholomorphic mappings by automorphisms of $\C^n$. 
Invent.\ Math.\ {\bf 112}, 323-349 (1993).
Correction, Invent.\ Math.\ {\bf 118}, 573-574 (1994).

\ii{[Gl1]} J.\ Globevnik: 
A bounded domain in $\C^N$ which embeds holomorphically into
$\C^{N+1}$. 
Ark.\ Mat.\ {\bf 35}, 313--325 (1997).

\ii{[Gl2]} J.\ Globevnik:  On Fatou-Bieberbach domains. 
Math.\ Z.\ {\bf 229}, 91--106, (1998).

\ii{[Gl3]} J.\ Globevnik:  
Interpolation by proper holomorphic embeddings of the disc into $\C^2$.
Math.\ Res.\ Lett {\bf 9}, 567--577 (2002).

\ii{[GS]} J.\ Globevnik, B.\ Stens\o nes:
Holomorphic embeddings of planar domains into $\C^2$.
Math.\ Ann.\ {\bf 303}, 579--597 (1995). 

\ii{[GSt]} R.\ E.\ Gompf, A.\ I.\ Stipsicz: $4$-manifolds and Kirby Calculus. 
American Mathematical Society, Providence,  1999.

\ii{[Gra1]} H.\ Grauert:
Holomorphe Funktionen mit Werten in komplexen Lieschen Gruppen.
Math.\ Ann.\ {\bf 133}, 450--472 (1957).

\ii{[Gra2]} H.\ Grauert: 
Analytische Faserungen \"uber holomorph-vollst\"andigen R\"aumen.
Math.\ Ann.\ {\bf 135}, 263--273 (1958).

\ii{[GRe]} H.\ Grauert, R.\ Remmert:
Theory of Stein Spaces.
Grundl.\ Math.\ Wiss.\ {\bf 227}, Springer, New York, 1977.

\ii{[Gn]} R.\ Greene: Isometric embeddings of Riemannian and
pseudo-Rie\-man\-nian manifolds.
Mem.\ Amer.\ Math.\ Soc.\ {\bf 97}, Providence, 1970.

\ii{[GW]} R.\ Greene, H.\ Wu:
Whitney's imbedding theorem by solutions of elliptic equations
and geometric consequences.
Proc.\ Symp.\ Pure Math.\ XXVII, 2, 287--297.
Amer.\ Math.\ Soc., Providence, 1975. 

\ii{[GSh]} R.\ Greene, K.\ Shiohama:
Diffeomorphisms and volume preserving embeddings
of non-compact manifolds.
Trans.\ Amer.\ Math.\ Soc.\ {\bf 255}, 403--414 (1979).

\ii{[Gro1]} M.\ Gromov: Convex integration of differential relations, I.  
Izv.\ Akad.\ Nauk SSSR {\bf 37}, 329--343 (1973). 
English transl.: Math.\ USSR Izv.\ 7 (1973).

\ii{[Gro2]} M.\ Gromov: Partial differential relations.
Ergebnisse der Mathematik und ihrer Grenzgebiete (3) {\bf 9},
Springer, Berlin--New York, 1986.

\ii{[Gro3]} M.\ Gromov:
Oka's principle for holomorphic sections of elliptic bundles.
J.\ Amer.\ Math.\ Soc.\ {\bf 2}, 851-897 (1989).

\ii{[GE]} M.\ Gromov, Y.\ Eliashberg: 
Nonsingular maps of Stein manifolds.
Func.\ Anal.\ and Appl.\ {\bf 5}, 82--83 (1971).

\ii{[GN]} R.\ C.\ Gunning, R.\ Narasimhan: 
Immersion of open Riemann surfaces. 
Math.\ Ann.\ {\bf 174}, 103--108  (1967). 

\ii{[GR]} C.\ Gunning, H.\ Rossi:
Analytic functions of several complex variables.
Prentice--Hall, Englewood Cliffs, 1965.

\ii{[Hae]} A.\ Haefliger: Vari\'et\'es feuillet\'es. 
Ann.\ Scuola Norm.\ Sup.\ Pisa (3) {\bf 16},  367--397 (1962). 

\ii{[Ham]} R.\ Hamilton: The inverse function theorem of Nash and Moser.
Bull.\ Amer.\ Math.\ Soc.\ {\bf 7}, 65--222 (1982).

\ii{[HL]} P.\ Heinzner, F.\ Kutzschebauch:
An equivariant version of Grauert's Oka principle.
Invent.\ Math.\ {\bf 119}, 317--346 (1995).

\ii{[HL1]}  G.\ Henkin, J.\ Leiterer:
Proof of Oka-Grauert principle without the induction over basis dimension.
Preprint, Karl Weierstrass Institut f\"ur Mathematik,
Berlin, 1986.

\ii{[HL2]} G.\ Henkin, J.\ Leiterer:
The Oka-Grauert principle without induction over the basis dimension.
Math.\ Ann.\ {\bf 311}, 71--93 (1998).

\ii{[Hi1]} M.\ Hirsch: Immersions of manifolds.
Trans.\ Amer.\ Math.\ Soc.\ {\bf 93}, 242--276 (1959).

\ii{[Hi2]} M.\ Hirsch: 
On embedding differential manifolds into Euclidean space.
Ann.\ Math.\ {\bf 73}, 566--571 (1961).

\ii{[H\"o]} L.\ H\"ormander:
An Introduction to Complex Analysis in Several Variables, 3rd ed.
North Holland, Amsterdam, 1990.

\ii{[Hu]} D.\ Husemoller: Fibre bundles. Third edition.
Graduate Texts in Mathematics, 20. Springer-Verlag,
New York, 1994.

\ii{[L\'ar]} F.\ L\'arusson:
Excision for simplicial sheaves on the Stein site
and Gromov's Oka principle. Preprint, December 2000.

\ii{[Lau]} H.\ B.\ Laufer:
Imbedding annuli in $\bf C^2$. 
J.\ Analyse Math.\ {\bf  26}, 187--215 (1973).

\ii{[Lei]} J.\ Leiterer:
Holomorphic Vector Bundles and the Oka-Grauert Priciple.
Encyclopedia of Mathematical Sciences {\bf 10}, 63--103,
Several Complex Variables IV, Springer, Berlin-New York, 1989.

\ii{[LV]} J.\ Leiterer, V.\ V\^aj\^aitu:
A relative Oka-Grauert principle on $1$-convex spaces.
Preprint, 2002.

\ii{[Na1]} R.\ Narasimhan: Holomorphic mappings of complex spaces. 
Proc.\ Amer.\ Math.\ Soc.\ {\bf 11}, 800-804 (1960).

\ii{[Na2]} R.\ Narasimhan: 
Imbedding of holomorphically complete complex spa\-ces. 
Amer.\ J.\ Math.\ {\bf 82}, 917--934 (1960). 

\ii{[Na3]} R.\ Narasimhan: On imbedding of Stein manifolds. 
Publ.\ Ramanujan Inst.\ {\bf 1}, 155-157 (1968/1969).

\ii{[N1]} J.\ Nash: $\cC^1$-isometric embeddings.
Ann.\ Math.\ {\bf 60}, 383--396 (1954).

\ii{[N2]} J.\ Nash: The imbedding problem for Riemannian manifolds.
Ann.\ Math.\ {\bf 63}, 20--63 (1956).

\ii{[Oka]} K.\ Oka: Sur les fonctions des plusieurs variables. III:
Deuxi\`eme probl\`eme de Cousin.
J.\ Sc.\ Hiroshima Univ.\ {\bf 9}, 7--19 (1939).

\ii{[Ph1]} A.\ Phillips:
Submersions of open manifolds.
Topology {\bf 6}, 170--206 (1967).

\ii{[Ph2]} A.\ Phillips: Foliations on open manifolds.
Comm.\ Math.\ Helv.\ {\bf 44}, 367--370 (1969).

\ii{[Ph3]}  A.\ Phillips: Maps of constant rank. 
Bull.\ Ann.\ Math.\ Soc.\ {\bf 80}, 513--517 (1974).

\ii{[Pr1]} J.\ Prezelj: 
Interpolation of embeddings of Stein manifolds on discrete sets. 
Math.\ Ann., to appear. 

\ii{[Pr2]} J.\ Prezelj:	Weakly regular embeddings of Stein spaces
with isolated singularities.
Preprint, 2002.

\ii{[Ram]} K.\ J.\ Ramspott: \"Uber die Homotopieklassen holomorpher 
Abbildungen in homogene komplexe Mannigfaltigkeiten. 
Bayer.\ Akad.\ Wiss.\ Math.-Natur.\ Kl.\ S.-B.\ 1962, 
Abt.\ II, 57--62 (1963).

\ii{[Rem]} R.\ Remmert:
Sur les espaces analytiques holomorphiquement s\'epa\-rab\-les et
holomorphiquement convexes. 
C.\ R.\ Acad.\ Sci.\ Paris {\bf 243}, 118--121 (1956). 

\ii{[RR]} J.-P.\ Rosay, W.\ Rudin:
Holomorphic maps from $\C^n$ to $\C^n$.
Trans.\ Amer.\ Math.\ Soc.\ {\bf 310}, 47--86 (1988).

\ii{[Sht]} U.\ Schaft: Einbettungen Steinscher Mannigfaltigkeiten.
Manuscripta Math.\ {\bf 47}, 175--186 (1984).

\ii{[Sch]} J.\ Sch\"urmann:
Embeddings of Stein spaces into affine spaces of minimal dimension.
Math.\ Ann.\ {\bf 307}, 381--399 (1997).

\ii{[Sl1]} Z.\ Slodkowski: Holomorphic motions and polynomial hulls. 
Proc.\ Amer.\ Math.\ Soc.\ {\bf 111}, 347--355 (1991). 

\ii{[Sl2]} Z.\ Slodkowski: Extensions of holomorphic motions. 
Ann.\ Scuola Norm.\ Sup.\ Pisa Cl.\ Sci.\ (4) {\bf 22}, 185--210 (1995). 

\ii{[Sm1]} S.\ Smale: 
A classification of immersions of the two-sphere. 
Trans.\ Amer.\ Math.\ Soc.\ {\bf 90}, 281--290 (1958).

\ii{[Sm2]} S.\ Smale; The classification of immersions of spheres 
in Euclidean spaces.
Ann.\ of Math.\ (2) {\bf 69}, 327--344 (1959).
                                     
\ii{[Sp]} D.\ Spring: Convex integration theory 
(Solutions to the $h$-principle in geometry and topology).
Monographs in Mathematics {\bf 92}, Birk\-h\"auser, Basel, 1998.

\ii{[Ste]} J.-L.\ Stehl\'e: Plongements du disque dans $\C^2$.
(S\'eminaire P.\ Lelong (Analyse), pp.\ 119--130), 
Lect.\ Notes in Math.\ {\bf 275}, Springer, Berlin--New York, 1970. 

\ii{[Stn]} K.\ Stein: 
Analytische Funktionen mehrerer komplexer
Ver\"anderlichen zu vorge\-ge\-benen Periodizit\"atsmoduln 
und das zweite Cousinsche Pro\-blem.
Math.\ Ann.\ {\bf 123}, 201--222 (1951).

\ii{[Th1]} W.\ Thurston: 
The theory of foliations of codimension greater than one. 
Comm.\ Math.\ Helv.\ {\bf 49}, 214--231 (1974).

\ii{[Th2]} W.\ Thurston: Existence of codimension one foliations.
Ann.\ Math.\ {\bf 104}, 249--268 (1976).

\ii{[Wi]} J.\ Winkelmann:
The Oka-principle for mappings between Riemann surfaces.
L'Enseignement Math.\ {\bf 39}, 143--151 (1993).

\bigskip
\ni Institute of Mathematics, Physics and Mechanics,
University of Ljubljana, Jadranska 19, 1000 Ljubljana, Slovenia

\smallskip
\ni E-mail: franc.forstneric@fmf.uni-lj.si

\bye